\documentclass[12pt,a4paper,twoside]{amsart}

\usepackage[T1]{fontenc}
\usepackage{amsmath}
\usepackage{amssymb}

\usepackage{latexsym} 

\newcommand{\newstuff}{} 
\newcommand{\newtext}{}

\newcommand{\newwtext}{} 
\newcommand{\newwtextt}{} 
\newcommand{\newtekst}{}

\newcommand{\clif}{\hbox{Clif}\,}
\newcommand{\R}{{\mathbb R}} 

\newcommand{\C}{{\mathbb C}} 
 
\newcommand{\cal}{\mathcal } 
\newcommand{\proofbox}{\hfill{$\Box$}}

\newcommand{\ud}{{\underline{\delta}}}

\newcommand{\Cnull}{{{\stackrel{\circ}C}{}^\infty}} 
\def\tilde{\widetilde}
\def \bfo {\begin {eqnarray*} }
\def \efo {\end {eqnarray*} }
\def \ba {\begin {eqnarray*} }
\def \ea {\end {eqnarray*} }
\def \beq {\begin {eqnarray}}
\def \eeq {\end {eqnarray}}
\def \supp {\hbox{supp}\,}
\def \diam {\hbox{diam }}
\def \rad {\,\hbox{rad}\,}
\def \ind {\hbox{Ind}\,}
\def \dist {\hbox{dist}\,}

\def\bra{\langle}
\def\cet{\rangle}
\def\bbra{\langle\langle}
\def\ccet{\rangle\rangle}
\def \e {\varepsilon}
\def \p {\partial}
\def \a {\alpha}
 
\def\la{\lambda}
\def\La{\Lambda}

\def\F{{\mathcal F}}

\def\hf{{\hat f}}
\def\hh{{\hat h}}
\def\ud{{\underline \delta}}
\def\Z{{\Bbb Z}}
\newtheorem{definition}{Definition}[section] 
\newtheorem{theorem}[definition]{Theorem} 
\newtheorem{lemma}[definition]{Lemma} 
 
\newtheorem{proposition}[definition]{Proposition} 
\newtheorem{corollary}[definition]{Corollary} 
\newtheorem{remark}[definition]{Remark}

\begin{document}
\title
[Inverse  Problems for Dirac operators]
{Inverse  Problems and Index Formulae for Dirac Operators} 
\author{Yaroslav Kurylev}
\address{
Yaroslav Kurylev,
Loughborough University,
Department of Mathematical Sciences, Leicestershire LE11 3TU, UK}
\email{}
\author{Matti Lassas}
\address{Matti Lassas, Helsinki University of Technology,
Institute of Mathematics,
PO Box 1100, FIN--02015 TKK, Finland}
\email{}

%%%%%%%%% To prevent amsart from capitalising the authors %%%%%%%%%%%%%%%%%%

%\makeatletter
%\def\@setauthors{%
%  \begingroup
%  \trivlist
%  \centering \@topsep30\p@\relax
%  \advance\@topsep by -\baselineskip
%  \item\relax
%  \andify\authors
%  \def\\{\protect\linebreak}%
%{\authors}%
%  \endtrivlist
%  \endgroup
%}
%\makeatother

%%%%%%%%%%%%%%%%%%%%%%%%%%%%%%%%%%%%%%%%%%%%%%%%%%%%%%%%%%%%%%%%%%%%%%%%%%%

\maketitle

%
%
%\begin{document} 
%\title{Inverse  Problems and Index Formulae for Dirac Operators} 
%
%\author{$\begin{array}{c} \hbox{Yaroslav Kurylev}
%\\  \hbox{University of Loughborough}\\ \hbox{UK} \end{array}$\quad
%$\begin{array}{c} \hbox{Matti Lassas}
%\\ \hbox{ Helsinki Univ.\ of Technology}\\ \hbox{Finland} \end{array}$
%} 
%\date{}%\today}
%\maketitle 

\noindent {\bf  Abstract:} 
{\it 
We consider a
Dirac-type operator $D_P$  on a vector bundle  $V$ 
over a compact Riemannian  manifold $(M, g)$ with a nonempty boundary.
The operator $D_P$ is specified by a boundary condition $P(u|_{\p M})=0$
where  $P$ is a projector which may be a non-local, i.e.\ a
 pseudodifferential operator.
We assume the existence of 
a chirality operator 
{\newwtextt which decomposes $L^2(M, V)$ into two orthogonal
subspaces
$X_+  \oplus X_-$.
Under certain conditions,  the operator $D_P$ 
restricted to $X_+$ and $ X_-$ defines  a pair of 
Fredholm operators which maps $X_+\to
 X_-$ and  $X_-\to  X_+$ correspondingly,  giving rise to a superstructure on $V$. }
 In this paper we consider the questions of 
 determining the
 index of $D_P$ and  the reconstruction of $(M, g), \,V$ and $D_P$ from the 
 boundary data on $\p M$.
 The data used is either the Cauchy data, i.e. the restrictions to 
 $\p M \times \R_+$ of the solutions to the hyperbolic Dirac equation,
 or the boundary spectral data, i.e. the set of the eigenvalues and the boundary values of 
the eigenfunctions of $D_P$. We obtain  formulae for the index and 
prove uniqueness results
for the inverse boundary value problems. We apply the 
obtained results to the classical Dirac-type operator in $M\times \C^4, \, M \subset \R^3$.}


\smallskip
 

\section{Introduction}\label{introduction 1}

Recent study of inverse problems has successfully 
shown that
external measurements can be used to uniquely determine coefficients
of various partial differential equations  modeling macroscopic and
microscopic phenomena.  Examples of these are
the paradigm problem, the inverse problem for the conductivity 
equation encountered in the impedance tomography, the inverse problem
for the wave equation with an inhomogeneous wave speed,
 the inverse scattering
problem for the Schr\"odinger operator 
(for different approaches see e.g. \cite{Ast1,Be1,BeKu,IY,KaKuSo,Na1,NaSyU,SyU}. In
these problems  the structure of the  space is
{\it a priori} known before measurements, being usually a known
domain of an Euclidean space. Recently, many inverse problems 
have
 been generalized
to the  cases where the underlying space is not {\it a priori} known 
but is assumed
to be some unknown smooth Riemannian manifold. This was followed 
by the  change of the attitude from proving the uniqueness in the considered
inverse problem to describing the corresponding groups
of transformations of the unknown object which preserve the measured 
data, e.g.\ isometries of Riemannian manifolds,  gauge equivalences of 
differential operators, etc.  From this point of view any anisotropic inverse 
problem, even in a given domain $\Omega\subset \R^n$,
 falls into this second category
having a non-trivial group of transformations which consists of some boundary 
preserving diffeomorphisms of the domain, see
 e.g.\ \cite{Ast2,BeKu3,Ku1,KL1,LU,PU} and, for a detailed description 
of this approach, \cite{KKL}.

A  natural task is to generalize this approach further to
inverse problems on
general vector bundles as, 
in the  modern physics, physical phenomena are often modeled
by equations on  bundles. Indeed, the vector bundles
 are encountered  
from the early development of the relativistic quantum mechanics
to the modern quantum field theory. 
Thus, in this paper our aim is to bring the mathematical study of 
inverse 
problems  closer to the modern physics 
and to study  inverse problems for one of its  basic equations, 
the Dirac equation on a general vector bundle. 
Assuming the existence of a chirality operator,
 we consider a Dirac equation on an 
unknown bundle and reconstruct the underlying manifold, the bundle
structure on it and also the  coefficients of the equation 
up to a natural group of transformations, namely, the isometries
of the manifold and bundlemorphisms of  the bundle.
 To this end  we further advance
the techniques
developed for the inverse problems for Maxwell's equation 
\cite{KLS1,KLS2,KLS3}.

For previous results in inverse
 problems for the Dirac equation in the 1-dimensional case,
see e.g.\  \cite{DK,Ger,Horv} and in 
the 3-dimensional case, see 
\cite{Hac,Iso,Jung,Ts}. 
In the $3$-dimensional case all these papers deal with  perturbations
of the canonical Dirac operator $D_0$ in
 $\R^3$ 
or $\Omega\subset \R^3$,
\beq 
\label{euclid}
D_0u =i \sum_{k=1}^3 \alpha_k \p_k u + mc^2\alpha_0 u, \quad
u \in L^2(\R^3; \C^4),
\eeq 
where $\alpha_{\nu}, \nu=0,\dots, 3,$ are the standard Dirac
matrices. Thus the structure of the underlying space 
and Dirac bundle was known {\it a priori}.

%It should be noted that even though we study
%inverse boundary value problem, we hope that the results
%will be useful for the inverse scattering problems.

\smallskip

\noindent {\bf  1.1.\ Formulation of the problem.}
In this paper we consider two types of inverse boundary value
 problems for the Dirac-type operator, $D$, acting on sections
of a complex vector bundle, $V$, over a compact, connected
Riemannian $n$-manifold $(M,g)$ with a non-empty boundary,  $\p M$.
Roughly speaking  (for a rigorous formulation of the problems see
Subsection 2.4  below), the first inverse problem deals with waves, i.e.
the solutions to the hyperbolic problem,
\beq
\label{dyn}
(i\p_t+D)u(x,t)=0\quad\hbox{in }M\times \R_+  ,
\quad u|_{t=0}=0.
\eeq
Then the data used is the set 
of  the Cauchy data of these waves restricted to 
 $\p M \times\R_+$, i.e. the set 
%\beq
%\label{Cauchy}
$C_0(D)=\{u|_{\p M \times\R_+}:\ 
u\ \hbox{satisfies } (\ref{dyn})\}.$
%\eeq 
The corresponding inverse problem, called 
the {\it dynamic inverse problem}, 
is that of the unique determination of    the manifold $M$, the bundle $V$, and the operator $D$
when we are given the boundary $\p M$, the restriction $W$ of the bundle to $\p M$, 
 $W=V|_{\p M}$,
and the {\it dynamic boundary data} $C_0(D)$.

The second inverse problem or, more rigorously, a family of inverse 
problems, deals with  the self-adjoint
 Dirac-type operator, 
\beq
\label{P}
D_P u = Du; \quad 
{\mathcal D}(D_P)= \{u \in H^1(M, V): \, P(u|_{\p M}) =0\}.
\eeq
Here  $P$ is an orthoprojector in $L^2(\p M, W),$ 
given by a zero-order 
classical
pseudo-differential operator on $\p M$ which, in general, 
is not a local operator. 
The corresponding data are the {\it boundary spectral data}, 
\ba
%\label{BSD before}
(\lambda_k,\, \phi_k|_{\p M}),\quad k=1,2,\dots,
\ea
where $\lambda_k,\, \phi_k$ are the eigenvalues and normalized
eigenfunctions of $D_P$.

Throughout the paper, we assume the  existence of a chirality operator
$F:V\to V$, acting fiber-wise, i.e.,
$
F: \pi^{-1}(x) \to \pi^{-1}(x), 
$
with $\pi: V \to M$ being the projector onto the base $M$.
The chirality operator satisfies the anti-commutation relation
\beq
\label{F-comm}
D \circ F+F\circ D=0
\eeq
% (for further properties of $F$ see Subsection 2.1 ).
%Then we show that the Cauchy data of waves, i.e., the set (\ref{Cauchy}),
%or  Boundary spectral data for 
%the Dirac-type operator $D_P$ defined by 
%any appropriate projector $P$ in (\ref{P}), 
%determine the Riemannian manifold $(M,g)$, up to an isometry,
%and also the Dirac bundle $V$ and Dirac-type operator $D_P$, up
%to a bundlemorphism.
and makes it possible, under some conditions,
to split $D_P$ into a pair of operators $D_P^{\pm}$,
$$
D_P \hspace{-1mm}= \hspace{-1mm} \left( \hspace{-2mm}\begin{array}{cc} 
0 &D_P^-   \\ 
D_P^+ & 0
\end{array} \hspace{-2mm}\right).
$$
We will 
derive formulae to express 
$\ind \,(D_P^{\pm})$ in terms of the dynamical boundary data.

The main tool in our study of the inverse problems
for the Dirac-type equations is the boundary control  method,
see \cite{Be1}, or,  more precisely,
its geometric version, see \cite{KKL} for a detailed exposition. 
The method is based on the finite velocity of the wave propagation
and control properties of the waves generated by boundary sources.
% together
%with control properties, in the domains of influence of open subsets 
%of the boundary, of the waves generated by boundary sources.

%The non-local nature of the projectors $P$ seems to prevent from using
%the Boundary Control method to solve previously formulated
%inverse problems. In this paper, we use a 
%special Dirac-type operator 
%of form (\ref{P}) with
%{\it local boundary condition} to resolve this difficulty.  The corresponding
%projector, $P_{\Gamma}$, is of the form,
%\beq
%\label{PGamma}
%P_{\Gamma} = \frac12 (I+\Gamma),
%\eeq
%where $\Gamma= \gamma(N) \circ F$, 
%with  $\gamma(N)$ being  a skew-adjoint operator acting on the sections of
%$V|_{\p M}$ (see section 2.1 for a detailed explanation).

The plan of the paper is as follows:
In Section 2 we give definitions of the Dirac
bundles and Dirac-type operators,
introduce a decomposition of
$L^2(M, V)$ into two channels $X_+ \oplus X_-$ 
and formulate rigorously
the main results of the paper. In Section 3 we describe properties
of the hyperbolic Dirac-type equation 
(\ref{dyn}) and provide a unique continuation result. 
In Section 4 we derive 
a formula  for the energy of a wave in each
channel in terms of the boundary data 
and prove some global and local boundary controllability
results for the hyperbolic Dirac-type equation. 
Generalized boundary sources and an extension of the
 time-derivative $\p_t$ 
are introduced in Section 5. The relation of  $\p_t$ and $D_P$
is studied in Section 6 where we obtain 
index formulae for $D_P$ in terms of the boundary data. In Section 7 we 
prove that the dynamic inverse data determine the manifold $M$
and its metric $g$. In Section 8 we introduce focusing sources which
generate waves that, at given time, are localized
at a single point $x\in M$. These waves  with
different polarizations are used in the reconstruction of the bundle $V$
 and operator $D_P$
 from the dynamic inverse data. 
In Section 9 we study the 
high-energy asymptotics  of the solutions to prove the uniqueness for
the inverse boundary spectral problem. 
In Section 10 we apply the obtained results
to the Dirac-type equation (\ref{euclid}) in
$\R^3$  leading
to an analog of the Aharonov-Bohm effect. 
%{\newtext 
%In appendix A we show that a special boundary projector $P_{\Gamma}$, called
%a {\it local boundary projector},
% gives rise to a self-adjoint Dirac-type operator. At last, Appendix B contains a sketch of the proof
% of the finite velocity propagation for the hyperbolic Dirac operator 
% (\ref{dyn}) with this local boundary condition.
%}
 In Appendix we prove, using the construction of Br\"uning and Lesch,
the self-adjointness of the Dirac-type operator with
the local boundary condition.
%and consider the finite velocity of the wave propagation.

\section{Definitions and main results.}
\label{dirac bundle}

\noindent {\bf  2.1.\ Definitions.} 
Here we introduce  some basic definitions, results and 
examples about 
Dirac bundles and operators on them
following the terminology of \cite{BrL,FS,Taylor}.

Let $(M,g)$ be a compact connected $C^\infty$-smooth Riemannian $n$-mani\-fold with a non-empty boundary.
Let $V$ be a smooth complex vector bundle over $M$
and denote the projection onto the base manifold by $\pi:V\to M$. 
Each fiber $\pi^{-1}(x)$ is a complex $d$-dimensional vector
space with a Hermitian inner product
$\bra\,\cdot ,\,\cdot \cet=\bra\,\cdot ,\,\cdot \cet_x$. 
We define $|\phi|_{x}^2=\bra \phi,\phi
\cet_{x}$ and denote  the smooth sections of $V$ by $C^\infty(M,V)$.
Endomorphisms $\hbox{End}\,(V)$ are fiber-preserving smooth maps,
$L:V\to V$, $\pi L(\phi)=\pi \phi$ that are  linear in each fiber.
By $\hbox{Clif}(M)$ we denote the Clifford bundle over $M$.
This means that, at each $x\in M$, the fiber $\hbox{Clif}_x(M)$ 
is an algebra generated by the vectors in $T_xM$ with a product $\cdotp$ 
satisfying the relation
\beq
\label{clifford}
v\,\cdotp w+w\,\cdotp v=-2g(v,w) \quad \hbox{for}\ v,w\in T_xM\subset 
\hbox{Clif}_x(M).
\eeq
We assume that there is a fiber-wise  map 
$\gamma:\hbox{Clif}(M)\to \hbox{End}\,(V)$
that provides  a Clifford module structure for  $V$, i.e.\ 
$\gamma$ gives an action, for any $x \in M$, 
of the algebra $\clif_x(M)$ on fibers 
$\pi^{-1}(x)$ of $V$ and
satisfies
\beq\label{commutator of gamma}
\bra \gamma(v)\lambda, \mu \cet_x + \bra \lambda, \gamma(v)\mu \cet_x=0,
\quad  v \in T_x(M), \,\, \lambda,\mu \in \pi^{-1}(x).
\eeq
Let
$\nabla$ be the Levi-Civit\'a connection in $(M,g)$. 
It defines a connection on $\clif(M)$ that satisfies
\ba
\nabla_X(v\,\cdotp w)=(\nabla_X v)\cdotp w+v\,\cdotp (\nabla_Xw),
\ea
where $X\in C^\infty(M,TM)$ and $v, w\in C^\infty(M,\clif(M)).$
We denote by  $\nabla$ also a connection on $V$ 
that is compatible with $\bra\cdotp,\cdotp\cet$ and $\gamma$,
i.e.,
\ba
& &X\bra \phi,\psi\cet =\bra \nabla_X\phi,\psi\cet+
\bra \phi,\nabla_X\psi\cet,
\\
& &\nabla_X(\gamma(w)\phi)=\gamma(\nabla_X w)\phi+\gamma( w)\nabla_X\phi,
\ea
where $X\in C^\infty(M,TM),$ $w\in C^\infty(M,\clif(M)),$ and
$\phi,\psi\in C^\infty(M,V).$
When these conditions are satisfied, we say that
$(V,\bra\,\cdot ,\,\cdot \cet,\gamma,\nabla)$ is a Dirac bundle 
over $(M,g)$, 
{\newtext for brevity, a Dirac bundle $V$.}

The unperturbed Dirac operator, $D_0$ on the Dirac bundle $V$,
which is sometimes called the {\it Dirac operator associated with a 
Clifford connection},  is 
locally given by 
\beq\label{dirac 1}
D_0 u|_U= \sum_{j=1}^n \gamma(e_j)\nabla_{e_j}u|_U,
\eeq
where $(e_1,\dots,e_n)$ is 
 an orthonormal frame in an open set $U\subset M$. 

As noted before, we assume that  there is
a chirality operator $F:V\to V$ acting fiber-wise and satisfying
\beq\label{chirality2}
 & &\bra F\lambda,\mu\cet_x=\bra \lambda,F\mu\cet_x, 
\quad  
\,\, \lambda,\mu \in \pi^{-1}(x),
\\ \nonumber
& &\gamma(v)\circ F+F\circ \gamma(v)=0, \quad v \in T_xM, 
\quad F^2=I,
\\ \nonumber
& &\nabla F=0,\quad \hbox{i.e.}\quad \nabla_X(Fu) =  F(\nabla_X u)
\eeq
where $
X\in C^\infty(M,TM)$ and $ u \in C^\infty(M,V).$
When  $n=\hbox{dim}(M)$ is even, a chirality operator always
exists, namely,
\beq\label{ch op}
F\phi|_U=
(\sqrt{-1})^{n/2} \gamma(e_1)\gamma(e_2)\dots\gamma(e_n)\phi|_U.
\eeq
A chirality operator exists also when we deal with
a space-type hypersurface on a Lorentzian manifold, e.g.\  \cite{HMR}.

Using a chirality operator, we define the fiber-wise 
orthogonal projectors onto the $+1$ and $-1$ eigenspaces of $F$,
\beq
\label{30.11.4}
\Pi_+=\frac 12 (I+F),\quad\Pi_-=\frac 12 (I-F).
\eeq
We denote by  $Q\in \hbox{End}\,(V)$  a self-adjoint 
potential, i.e.,
\beq\label{Q self-adjoint}
\bra Q\lambda,\mu\cet_x=\bra \lambda,Q\mu\cet_x,\quad \lambda,\mu
\in \pi^{-1}(x),
\eeq
which respects the chirality structure, i.e.\
\beq\label{respect}
F \circ Q+Q\circ F=0.
\eeq
We call the minimal Dirac-type operator, $D_{\rm min}$, the operator
\beq
\label{minimal2}
D_{\rm min}u =\, D_0 u + Qu, \quad u \in {\mathcal D}(D_{\rm min})
=C_0^{\infty}(M, \,V),
\eeq
where $ C_0^{\infty}(M, \,V)$ is the space of the smooth
compactly supported section of $V$. 
%We note that this is a special class
%of the operators considered in \cite{BrL}.

\smallskip

\noindent {\bf  Example 1.}
 As an example, we consider the
 so-called form-Dirac operator, see e.g. \cite{FS,Taylor}.
 Let $(M,g)$ be an oriented Riemannian 
manifold of  dimension $n$ and 
$\Omega^j M$ be the space of the 
complex-valued differential  $j$-forms on $M$. Let
$\Omega M=\Omega^0 M\oplus
\dots\oplus \Omega^n M$ be the Grassmannian bundle on $M$
and $*$ the Hodge operator of $\Omega(M)$ related to the metric $g$.
For $\lambda\in \Omega_x^j M$ and
$\mu\in \Omega_x^kM$,  let $\bra \lambda,{\mu}\cet_x=
0$ if $j\not =k$ and  $\bra \lambda,\mu\cet_x=
*(\lambda \wedge *\overline \mu)$ for $j=k$.

Denote by $I:T_xM\to T_x^*M$ the identification 
$I(a^j\frac \p{\p x^j}) 
=g_{jk}a^jdx^k$ and by $\iota_{v}:\Omega^j_x M\to \Omega^{j-1}_x M$
 the  inner product with a vector $v \in T_xM$, e.g.
 \bfo
\iota_{v}w(v_1, \dots,v_{j-1})=w(v, v_1, \dots, v_{j-1}),
\efo
for any $w \in \Omega_x^jM, \, v, v_1, \dots, ,v_{j-1} \in T_xM$.
Then the map $\gamma:T_xM \to \hbox{End}(\Omega_x M)$,
\ba
\gamma(v) \lambda=-I(v)\wedge \lambda+\iota_{v} \lambda 
\ea 
extends to a homomorphism $\gamma: \hbox{Clif}_x(M)\to \hbox{End}
(\Omega_x M)$.
Let $\nabla$ be the Levi-Civita connection on $\Omega M$.
Then 
$(\Omega M,\bra\,\cdot ,\,\cdot \cet,\gamma,\nabla)$ is
a Dirac bundle. 
The unperturbed Dirac operator  (\ref{dirac 1})
on this bundle is 
\ba
%\label{form-Dirac}
\quad \quad D_0:\Omega M\to \Omega M, \quad
%\\  \nonumber
D_0=d+\delta: \,
\Omega^j M \, \rightarrow \, \Omega^{j+1}M\oplus  \Omega^{j-1}M,
\ea
where $d$ is the exterior differential and
\bfo
\delta=(-1)^{nj+1}*d*:\Omega^j M\to \Omega^{j-1}M
\efo 
is the 
codifferential.
A chirality operator, $F$ for the form-Dirac operator may be defined 
simply as
\beq\label{Form F}
F \lambda= (-1)^j \lambda, \quad \lambda \in \Omega^jM,
\eeq
splitting the Grassmannian bundle into the forms of the even and odd orders,
$
\Omega M = \Omega^{\tiny e}M \oplus  \Omega^{\tiny o}M.
$
\smallskip

\noindent {\bf  Example 2.}
Consider Maxwell's equations in $M\times \R$, $M\subset \R^3$,
\ba
%\label{Maxwell-Faraday 1}
& & {\rm curl}\,E(x,t) = - B_t(x,t),\quad D(x,t)=\epsilon(x) E(x,t), \\
%\quad{\rm div}\,B(x,t)=0,\\
\nonumber%\label{Maxwell--Ampere 1}
& & {\rm curl }\,H(x,t) =\phantom{-} D_t(x,t),\quad B(x,t)=\mu(x) H(x,t)
% \quad   {\rm div}\,D(x,t)=0
\ea
where $\epsilon(x)$ and $\mu(x)$ are positive definite matrix valued functions.
The velocity of the wave propagation is independent of  the wave polarization
if and only if  $\mu(x)=\alpha(x)^2 \mu(x)$ with some scalar function
$\alpha$. In this case the travel time is determined by the
Riemannian metric $g_{ij}=
 \alpha^{-2}{\rm det}(\epsilon)^{-1} 
\delta^{ik}\epsilon^j_k.$
As shown in \cite{KLS1,KLS3}, Maxwell's system may be extended to the
hyperbolic equation
\beq\label{Dirac-Maxwell}
i\p_t \omega+i(d-\delta_\alpha)\omega=0,\quad \omega(x,t)=
(\omega^0,\omega^1,\omega^2,\omega^3)\subset \Omega M,
\eeq
where $\delta_\alpha \omega^j =\alpha \delta(\alpha^{-1}\omega^j)$
and $\delta$ is the codifferential with respect to the travel
time metric $g$. Indeed, a solution of Maxwell's system
gives rize to a solution of the Dirac-type equation
(\ref{Dirac-Maxwell}) with 
$\omega^0=0,$ $\omega^1=E_jdx^j,$ $\omega^2=
B_3dx^1\wedge dx^2-B_2dx^1\wedge dx^3+B_1dx^2\wedge dx^3$,
and $\omega^3=0$. Defining the inner product 
\ba
\bra \omega^j,\lambda^k\cet_x=\alpha(x)^{-1}
*(\omega^j\wedge *
{\overline \lambda^k})\delta^{jk},\quad \hbox{for }
\omega^j\in \Omega^j_xM,\ 
\lambda^k\in \Omega^k_xM,\ea the operator $i(d-\delta_\alpha)$
is formally self-adjoint in $\Omega M$.
This operator defines a Dirac-type bundle
$(\Omega M,\bra\,\cdotp,\cdotp\cet_x,\tilde \gamma,\nabla)$
where $\tilde \gamma(v)\lambda=-i(I(v)\wedge \lambda+\iota_v\lambda)$
and $\nabla$ is  the Levi-Civit\'a connection in $(M,g)$. 
The chirality operator (\ref{Form F}) is a chirality operator
on this bundle, too. 

\smallskip

\noindent {\bf  2.2.\ Self-adjoint operators.}
We return now to the discussion of the boundary conditions
{\newtext  used to extend the minimal Dirac operator 
(\ref{minimal2}) to a self-adjoint operator }
{\newtekst in $L^2(M, V)$ with the inner product
\ba
\bbra \phi,\psi\ccet=\int_M \bra \phi(x),\psi(x)\cet_x \,dV_g(x),
\ea 
where $dV_g$ is the Riemannian volume  on $(M,g)$.}
Self-adjoint extensions of $D_{\rm min}$
defined by  {\it non-local} boundary conditions,
go back to the centennial work by Atiyah-Patodi-Singer  \cite{APS}. In this 
paper, we will work with the self-adjoint extensions of $D_{\rm min}$
extensively studied by Br\"uning and Lesch
\cite{BrL}, see also  \cite{FS,HSZ1}.
% We note, however, that the
%class of the Dirac-type operators considered in \cite{BrL} is wider 
%than the Dirac-type operators with chirality associated with a 
%Clifford connection considered in this paper.

% We start with an introduction of some 
% function spaces used in the paper.
%Let $dV_g$ be the Riemannian volume  on $(M,g)$.
%For smooth sections $\phi,\, \psi \in C^{\infty}(M,\,V)$  we define 
%the $L^2$-inner product
%\ba
%\bbra \phi,\psi\ccet=\int_M \bra \phi(x),\psi(x)\cet_x \,dV_g(x)
%\ea 
%and the norm $\|\phi\|_{L^2(M,V)}^2=\bbra \phi,\phi\ccet$. 
%The completion
%of  $C^\infty(M,V)$ with respect to this norm
%is denoted by $L^2(M,V)$.
{\newtekst In the following we denote by $H^s(M,V)$
the Sobolev spaces   of  sections of $V$ 
with components 
 in the Sobolev spaces $H^s(U)$
in local smooth trivializations $\Phi_U:
\pi^{-1}(U)\to U\times \C^d$, $U\subset M$ and
  by  $H^s_0(M,V)$ the closure
 of  $C^{\infty}_0(M,\,V)$  in  $H^s(M,V)$. }

{\newtext  Using the natural embedding $j:\p M\to M$, we
introduce the induced bundle $W=V|_{\p M}=j^*V$ on $\p M$.}

Next we  consider the possible extensions, in $L^2(M, V)$,
of the minimal Dirac-type operator (\ref{minimal2}).
We use $D$ to denote the maximal extension
of the operator $D_{\rm min}$ to $H^1(M, V)$,
\beq
\label{dirac-type, maximal}
Du = (D_0+Q) u, \quad {\cal D}(D) = H^1(M, V).
\eeq
It is shown in \cite{BrL} that
there is a 
wide class of self-adjoint extensions of $D_{\rm min}$
 defined by  boundary conditions of form (\ref{P}),
\beq
\label{dirac-type}
D_P u = D u, \quad {\cal D}(D_P) = \{u \in H^1(M, V): \, P (u|_{\p M})=0\},
\eeq
so that $D_P\subset D$.
Here $P$ is a  zero-order 
classical
pseudo-differential operator defining an orthoprojector
 on $L^2(\p M,W)$. It should satisfy the following conditions:
\beq
\label{anticommutation}
P \gamma(N) = \gamma(N)(I-P),
\eeq
and 
\beq
\label{Fredholm}
\{P,\, P_{APS}\} \,\,\hbox{is a Fredholm pair},
\eeq
that is, $ P_{APS}:\hbox{Ran}\,(P)\to \hbox{Ran}\,(P_{APS})$ is a 
Fredholm operator. Here $N$  stands for the unit interior normal field of $\p M$ and 
the operator $P_{APS}$  is associated with 
 the {\it Atiyah-Patodi-Singer} boundary condition.
Namely, let  $A(0)$, sometimes called the {\it hypersurface
Dirac operator}, (see e.g.\  Appendix or \cite{HSZ1}),
 be a self-adjoint operator
in $L^2(\p M, W)$, 
\beq
\label{A(0)}
A(0) = 
- \gamma(N) \sum_{\a=1}^{n-1} \gamma(e_{\a})\triangledown_{\a} +
\frac{n-1}{2} H(x'), \quad \triangledown_{\a}=\triangledown_{e_\a}.
\eeq
Here the vector fields $e_{\a}, \, \a=1, \dots, n-1,$ form a
local orthonormal frame
 on $\p M$, while $H(x')$ is the mean curvature
of $\p M$ at the point $x'$, where $x'=(x_1, \dots, x_{n-1})$
are local coordinates on $\p M$. Then $P_{APS}$  is a spectral projector
of $A(0)$ such that 
\beq
\label{APS}
& &P_{(0,\infty)} \subset P_{APS} \subset P_{[0,\infty)},\\
& &P_{APS} = -\gamma(N)(I-P_{APS})\gamma(N) ,
\eeq
where we denote by  $P_{(a,\,b)}$ the spectral projector of $A(0)$ associated 
with an interval $(a,b)$. 

\begin{theorem} [Br\"uning-Lesch]
 \label{self-adjointenss}
Let 
$P$ be a  zero-order 
classical
pseudo-differential operator defining an orthoprojector which
satisfies conditions  
(\ref{anticommutation}) and (\ref{Fredholm}).
Let $D_{P}$  be the Dirac-type operator (\ref{dirac-type}).
Then
$D_{P}$ is self-adjoint, the spectrum of $D_{P}$
is discrete, and all eigenspaces are finite
dimensional.
Moreover, $D_P$ is regular, that is,
$D_Pu\in H^s(M,V)$ implies that
$u\in H^{s+1}(M,V)$ for any $s\geq 0$.
\end{theorem}
The proof of Theorem \ref{self-adjointenss} is given in  
\cite[Thm. 1.5]{BrL}
for more general Dirac-type operators
than those considered in this paper  and in  \cite{FS,HSZ1} for more restricted cases. 

In the sequel we will denote the eigenvalues of $D_P$,
numerated according to the their multiplicity, 
by $\la_k^P,\,
|\la_{k}^P| \leq |\la_{k+1}^P|$, and the corresponding orthonormal 
eigenfunctions by $\phi_k^P$. When there is no danger of confusion
we will skip the index $P$ in these notations.

Using projections $\Pi_+$ and $\Pi_-$, see 
(\ref{30.11.4}) we define the bundles
$V_+=\Pi_+ V$ and $V_-=\Pi_- V$
over $M$ with the projections to $M$ denoted by $\pi_{\pm}$.
Note that, for any $x \in M$,
$\pi_+^{-1}(x)$ and  $\pi_-^{-1}(x)$ are orthogonal  so that
$V=V_+\oplus V_-$ and 
$L^2(M,V)=L^2(M,\Pi_+V)\oplus L^2(M,\Pi_-V)$. 
With a slight abuse of notation we use
$\Pi_\pm$ for the orthoprojectors in $L^2(M,V)$ onto
$L^2(M,\Pi_\pm V)$.
The maximal Dirac-type operator $D$
satisfies $DF+FD=0$ so
that $D$
can be decomposed 
as 
\beq\label{splitting}\quad 
%\\ \nonumber 
& &D: \hspace{-1mm} 
H^1(M,\Pi_+V)  
\oplus
H^1(M,\Pi_-V)\hspace{-1mm} 
\to \hspace{-1mm} 
L^2(M,\Pi_-V)
\oplus
L^2(M,\Pi_+V),\\
\nonumber& &
D \left(\begin{array}{c}  
u_+ \\ 
u_-
\end{array}\right)=
 \hspace{-1mm}= \hspace{-1mm} \left( \hspace{-2mm}\begin{array}{cc} 
0 &D   \\ 
D & 0
\end{array} \hspace{-2mm}\right)
 \left(\begin{array}{c}  
u_+ \\ 
u_-
\end{array}\right),
\quad u_{\pm} = \Pi_{\pm} u.
\eeq
When 
\beq\label{decomp 1.}\quad
\quad \Pi_+ \, {\mathcal D}(D_P) \subset {\mathcal D}(D_P)\quad\hbox{or, equivalently,}\ \
\Pi_-\, {\mathcal D}(D_P) \subset {\mathcal D}(D_P),
\eeq
decomposition (\ref{splitting}) gives rize to a
decomposition of $D_P$. Namely, for
$D_P^\pm =D_P:
\Pi_\pm\, {\mathcal D}(D_P)
 \to L^2(M,\Pi_\mp V),$ 
\beq\label{decomposition0}
& &D_P :  \Pi_+ 
{\mathcal D}(D_P) \oplus 
\Pi_- {\mathcal D}(D_P)\to
L^2(M,\Pi_+V)\oplus
L^2(M,\Pi_-V),
\\ \nonumber
& &D_P \left(\begin{array}{c}  
u_+ \\ 
u_-
\end{array} \hspace{-2mm}\right)
=
\hspace{-1mm} = \hspace{-1mm} \left(
 \hspace{-2mm}\begin{array}{cc} 
0 &D_P^-   \\ 
D_P^+ & 0
\end{array} \hspace{-2mm}\right) \hspace{-1mm}
 \left(\begin{array}{c}  
u_+ \\ 
u_-
\end{array}\right).
\eeq 
If (\ref{decomp 1.}) is valid and thus 
the decomposition (\ref{decomposition0}) is possible,
we say that the bundle $V=V_+\oplus V_-$ on $M$
has a {\it superstructure}.
It then follows from \cite{BrL} that  the resulting operators $D_P^{\pm}$ are Fredholm with,
in general, non-trivial indices.
For important relations of the index and the geometry of the bundle, see e.g.
\cite{APS,Bis,gilkey2,gilkey1}.

\smallskip

\noindent {\bf  2.3.\ Local boundary condition.}
As noted in Introduction, the boundary projection $P$ in (\ref{dirac-type})
is usually non-local. However,
for a Dirac-type operator with chirality, it is possible
to introduce
 a {\it local boundary 
condition} which makes it self-adjoint.
To this end, we consider
 a fiberwise operator $\Gamma$ on $W=V|_{\p M}$,
\beq
\label{Gamma}
\Gamma:W\to W,\quad \Gamma \lambda=
F\,\gamma(N)\, \lambda,
\quad \lambda \in \pi^{-1}(x),\quad x\in \p M.
\eeq
It is clear from properties (\ref{commutator of gamma}) and 
(\ref{chirality2})
 that 
\beq\label{prop of Gamma}\quad\quad
& &\Gamma^2=I,\quad \Gamma \circ F+F \circ \Gamma=0,
\\ \nonumber & &\bra \Gamma\lambda,\mu\cet_x=
 \bra\lambda,\Gamma \mu\cet_x, \,\, \lambda, \mu \in \pi^{-1}(x).
\eeq
Let 
\beq
\label{PGamma}
P_{\Gamma}=\frac 12 (I+\Gamma),
\eeq
be an orthoprojector acting fiber-wise on 
$\pi^{-1}(x),\, x \in \p M$.
Then $P_{\Gamma}$ satisfies conditions (\ref{Fredholm}) and 
(\ref{anticommutation}) and, therefore, defines a self-adjoint Dirac-type
operator
 $D_{\Gamma}:=D_{P_\Gamma}$. We note that this operator will play a crucial
role in the sequel. Because of this, we give in the Appendix 
a sketch of the proof of its self-adjointness.

%There is also another Dirac-type operator associated with $\Gamma$.
%It is defined by the projector,
%\bfo
%P= I-P_{\Gamma}.
%\efo
%It has the same properties as $D_{\Gamma}$ which may be seen
%from the fact that the projector $P_{\Gamma}$ goes into $(I-P_{\Gamma})$
%if we change $F$ into $-F$. In this paper, except for Section 7, it would
%not play an important role. 
%However, in the proof of Theorem \ref{main1}
%in Section 7, we will use it together with  $D_{\Gamma}$. In this case,
%to distinguish these two operators, we will denote  $D_{\Gamma}$ by
% $D_{\Gamma}^+$ and the operator associated with $(I-P_{\Gamma})$ by
% $D_{\Gamma}^-$. Similarly, we will use superscripts and subscripts $\pm$
%to distinguish various objects related to $D_{\Gamma}^+$ and $D_{\Gamma}^-$.
%For example, we will denote $P_{\Gamma}$ by $P_{\Gamma}^+$ 
%and $(I-P_{\Gamma})$ by $P_{\Gamma}^-$.
\smallskip

\noindent {\bf  Example 1, continued.}
Clearly, $P_{\Gamma}$ defines a local boundary condition for the 
form-Dirac operator $d+\delta$ determined in Example 1. In addition,
there are other local boundary condition closely related to the de Rham 
complexes over $(M, g)$.
Recall that $j: \p M \to M$ is a natural embedding.
The {\it relative boundary condition} is defined by
\beq
\label{relative}
P_r(u|_{\p M}): = j^*(u) =0,
\eeq
and the {\it absolute boundary condition} by
\beq
\label{absolute}
P_a(u|_{\p M}): = j^*(\iota_N(u)) =0.
\eeq
In the future we denote the corresponding self-adjoint form-Dirac operators
by $(d+\delta)_r$ and $(d+\delta)_a$.
\smallskip

\noindent {\bf  Example 3. }
Another important example of a Dirac operator on the Grassmannian bundle
$\Omega M$, for $n=2k$, is given by
\bfo
D_1u \,=\, (d+\delta)u,\quad u\in {\mathcal D}(D_1)=\{
u\in H^1(M,V):\ P_{APS}(u|_{\p M})=0\},
\efo
where $P_{APS}$ is
 the Atiyah-Patodi-Singer boundary condition defined in subsection 2.2,
see e.g.\  \cite{APS,gilkey1,Me1}.
This operator is called
the  signature operator. In this case the chirality
operator is given by (\ref{ch op}).

\smallskip

\noindent {\bf  2.4. Boundary data for inverse problems.} 
%We finish this section by introducing the inverse boundary data used in 
%this paper.
\begin{definition}
\label{BSD}
{\rm 
Consider
a Dirac bundle
$(V,\bra\,\cdot ,\,\cdot \cet,\gamma,\nabla)$
with chirality $F$.
The induced bundle structure
on $\p M$ is the collection
\beq
\label{b-structure}
\{\p M,\, g|_{\p M},\, W,\, \bra
\cdotp,\cdotp\cet_x|_{x\in \p M},\,  \gamma(N)|_{\p M},\, F|_{\p M}\}.
\eeq
Here $\p M$ is considered as a Riemannian manifold
with the differentiable structure induced by the embedding
$j:\p M\to M$ and the metric $g|_{\p M}=j^*g$.
The bundle
$W=V|_{\p M}$ 
is the induced bundle on $\p M$,
 $\bra\cdotp,\cdotp\cet_x$ the Hermitian structure
to $W$,
$\gamma(N)|_{\p M}$ and $F|_{\p M}$ are the restrictions
of the Clifford action $\gamma(N)$ and of
the chirality operator $F$ on $W$.
}
\end{definition}

\begin{definition}
\label{BSD 2}
{\rm Let $D_P$ be a self-adjoint  Dirac-type operator of form
(\ref{dirac-type}) on a Dirac bundle
$(V,\bra\,\cdot ,\,\cdot \cet,\gamma,\nabla)$
with chirality $F$.
The set
\beq
\label{BSD1}
\{(\lambda_k,\, \phi_k|_{\p M})\}_{k=1}^{\infty},
\eeq
where $\lambda_k$ and $ \phi_k$ are the eigenvalues and the orthonormal
eigenfunctions of $D_P$
is called the
boundary spectral data of $D_P$.}
\end{definition}

To define another type of boundary data, observe that
any Dirac-type operator (\ref{dirac-type}) is associated with a
hyperbolic
 initial-boundary value problem
\beq
\label{IBVP0}
& &(i\p_t+D)u(x,t) = 0\quad \hbox{in }M\times \R_+,
\\
\label{IBVP0line 2}
& &P (u|_{\p M \times \R_+}) = f \in P {\Cnull}(\p M\times \R_+, W)), 
\quad
u|_{t=0}=0,
\eeq
where ${\Cnull}(\p M\times \R_+,X)$ is the class of the smooth functions on
 $\p M\times \R_+$
(with values in $X$) which are equal to $0$ near $t=0$.
%and
%\ba
%P {\Cnull}(\R_+ \times\p M, W))={\Cnull}(\R_+, P(C^\infty(\p M,W)))
%\ea
%is the
% range of $P$,
We denote by $u=u^f$ the solution of (\ref{IBVP0})-(\ref{IBVP0line 2}).
Problem (\ref{IBVP0})-(\ref{IBVP0line 2})
 defines the response operator, $\La_P$,
\beq
\label{response}
\La_P f = u^f|_{\p M \times \R_+}.
\eeq
When $P=P_\Gamma$, we denote $\La_\Gamma=\La_{P_\Gamma}$.
Other data for the inverse boundary
problem could be the response operator $\La_P$.
However, we prefer to work with  data which is independent of $P$.
\begin{definition}
\label{Cauchydata}
{\rm Let $D$ be the maximal Dirac-type  
operator (\ref {dirac-type, maximal}).
Then the Cauchy data set of $D$ is the set 
\beq
\label{cauchy}\quad\quad
C_0(D)=\{u|_{\p M \times \R_+} : \,
u\in \Cnull( M \times \R_+,  V) \ \hbox{satisfies (\ref{dyn})}\}.
\eeq
}
\end{definition}

Later we prove the following  equivalence of different type of data.

\begin{lemma}\label{lem: equivalence of data} Assume that the 
induced bundle structure
(\ref{b-structure}) on $\p M$ is given. Then the
Cauchy data set $C_0(D)$ determines the map $\Lambda_\Gamma$ and vice versa. 
Moreover, the set $C_0(D)$ and 
a projector $P$ satisfying  
(\ref{anticommutation}) and (\ref{Fredholm})
determine
the map $\Lambda_P$ and vice versa.
\end{lemma}

\noindent {\bf  2.5.\ Main results.} First 
we formulate our main index formula.
\begin{theorem} 
\label{main0}
Let $V$ be a Dirac bundle
over $(M,g)$ with a chirality operator $F$. Let $D_P$ be a self-adjoint Dirac-type operator
of the form  (\ref{dirac-type}).

Assume that we are given
the induced bundle structure
(\ref{b-structure}) on $\p M$
and the response operator $\La_P$. Then these data determine,
in a constructive way, a Hilbert space
$\F$,
a pair of projectors $B_+$ and  $B_-$  in $\F$ with $B_+ \oplus B_-=I$,
and an unbounded self-adjoint operator
$\p_t: 
\F
\to  \F$ with  ${\mathcal D}(\p_t)\subset \F$ 
such that 

%\smallskip
%\noindent 1. There exists a isomorphism $S:X
%\to L^2(M,V)$ such that $S({\mathcal D}(\p_t))={\mathcal D}(D_P)$
%and $D_P=S(-i\p_t)S^{-1}$.

\smallskip

\noindent 1. The condition (\ref{decomp 1.}) for $D_P$ is valid
if and only if $B_+{\mathcal D}(\p_t)
\subset {\mathcal D} (\p_t)$.

\smallskip

\noindent 2. If decomposition (\ref{decomp 1.}) is valid
then the operator
%can be decomposed
\beq\label{d_t decompos.}
\p_t^+ : B_+{\mathcal D}(\p_t)\to B_-\F,\quad \p_t^+ u =\p_t u
%\p_t = \left(\begin{array}{cc} 
%0 &\p_t^-   \\ 
%\p_t^+ & 0
%\end{array}\right): B_+{\mathcal D}(\p_t)
 % \hspace{-1mm}
%\opus \hspace{-1mm}
% B_-{\mathcal D}(\p_t)\hspace{-1mm}
%\to \hspace{-1mm}
% B_+X
%\times\hspace{-1mm}
%B_-X
%\eeq
%where  the restrictions $\p_t^\pm$ of $\p_t$ are 
\eeq
is a Fredholm operator
and 
\beq\label{eq: Main index}
\ind \,(\p_t^+)=\ind\, (D_P^+).
\eeq
\end{theorem}

In Theorem \ref{main0} the Hilbert space $\F$ is
the completion of the functions $P C^\infty_0(\p M\times [0,T], W))$,
{\newwtextt where $T>0$ is sufficiently large. This completion is taken }
with respect to a seminorm $\|\,\cdotp\|_\F$ that is
explicitly determined by the Cauchy data set $C_0(D)$
(see (\ref{innerproduct}), (\ref{U-formul.}), (\ref{B slava}) below).
 The operator 
$\p_t$ in $\F$ is an extension of the time derivative
\ba
\p_t f=\frac {\p f} {\p t},\quad\hbox{for }f\in  P C^\infty_0(\p M\times [0,T], W)).
\ea 
%
%
%This theorem is proven in section 5, where we also  construct
% a map
%\ba
%J:  P C^\infty_0(\p M\times [0,T], W)) \to {\mathcal D}(\p_t)\subset \F
%\ea
%with a dense image, where $T>0$ is specified in Section 5 such that
%\beq\label{eq: Lboro}
%\p_t (J f)= J(f_t),\quad f\in P C^\infty_0(\p M\times [0,T], W)).
%\eeq
%It is equation (\ref{eq: Lboro}) which explains the notation $\p_t$ for
%the considered operator in $\F$. 

Next we formulate our main results on
the reconstruction of the manifold and the Dirac bundle.
To this end, let $V$ and $\tilde V$ be vector bundles
over diffeomorphic manifolds $M$ and $\tilde M$. 
Let $\ell:M\to \tilde M$ be a diffeomorphism  and $\pi:V\to M$ and   $\tilde \pi:\tilde V\to \tilde M$ 
be projections on the base spaces. A map $L:V\to \tilde V$
is called a {\it bundlemorphism} compatible with $\ell$ if
$\ell(\pi \phi)=\tilde \pi (L \phi),\quad \phi\in V$ and 
\ba
& &L:\pi^{-1}(x)\to \tilde \pi^{-1}(\ell(x))\hbox{ is a linear isomorphism for any}x\in M.
\ea
Dealing with Riemannian manifolds and Dirac bundles  on them
we need to extend the notation
of a bundlemorphism.

\begin{definition}\label{def: bundlem} 
{\rm Let
 $(V,\bra\, \cdot,\cdot\, \cet,\gamma,\nabla)$ and
  $(\tilde V,\bra\, \cdot,\cdot\, \cet,\tilde \gamma,\tilde \nabla)$
  be Dirac bundles over isometric Riemannian manifolds $(M,g)$
  and $(\tilde M,\tilde g)$ with chirality operators
  $F$ and $\tilde F$, and let $D$ and $\tilde D$ be
  Dirac-type operators of  form  (\ref{dirac-type, maximal}) on them.
  A bundlemorphism $L$ compatible with an isometry $\ell:M\to \tilde M$ 
  is called a {\it Dirac bundlemorphism} if
  \beq\label{eq: lboro A}
  & &\bra L\phi, L\psi \cet_{\ell(x)} = \bra \phi, \psi \cet_x, \quad \phi,\psi\in \pi^{-1}(x),\\
  \nonumber & & \tilde \gamma( d\ell(v))= L\gamma (v)L^{-1},\quad
  v\in T_xM,\\  \nonumber 
  & & \tilde F=L F L^{-1},\quad \tilde D=L D L^{-1}, 
\eeq
where $d\ell$ is the differential map of $\ell$.
}
\end{definition}

We turn now to some induced structures 
on the boundary.

\begin{definition}\label{boundary morphism}
{\rm Let 
 $(V,\bra\, \cdot,\cdot\, \cet,\gamma,\nabla)$ and
  $(\tilde V,\bra\, \cdot,\cdot\, \cet,\tilde \gamma,\tilde \nabla)$
be Dirac bundles over manifolds $M$ and $\tilde M$
with $W=V|_{\p M}$ and  $\tilde W=\tilde V|_{\p M}$. 
 Let $\p M$ and $\p {\tilde M}$ be isometric with 
an isometry
$\kappa:\p M\to \p\tilde M$. 
We say
that a bundlemorphism $K:W\to \tilde W$ is an 
{\it isomorphism between the
induced bundle structures} of $\p M$ and $\p \tilde M$ compatible
with  $\kappa$
if 
 \ba%\label{eq: lboro B}
%& &\kappa(\pi \phi)=\tilde \pi (K \phi),\quad \phi\in W,\\  \nonumber 
%& &K:\pi^{-1}(x)\to \tilde \pi^{-1}(\kappa(x))\hbox{ is an linear i%somorphism for }x\in \p M,\\
 %\nonumber 
  & &\bra K\phi, K\psi \cet_{\kappa(x)} = \bra \phi, \psi \cet_x,
  \quad \tilde F_{\kappa(x)}=K F_x K^{-1} ,
  \\
  \nonumber & & \tilde \gamma( d\kappa( v))
  = K\gamma (v)K^{-1},\quad
   \tilde \gamma(\tilde N(\kappa(x)))= K\gamma (N(x))K^{-1}.
 % \\ \nonumber
 % & & \tilde F=K F K^{-1} 
\ea
Here $x\in \p M,\, \kappa(x) \in \p {\tilde M}$,
$\phi,\psi\in \pi^{-1}(x)$ and $   v\in T_x(\p M).$ At last, 
$N(x)$ and $\tilde N(\tilde x)$ are the unit interior normal vectors to $\p M$ at $x$
and  to $\p \tilde M$ at $\tilde x$. 
}
\end{definition}

Let now $u\in  {\Cnull}(M\times \R_+, V)$ be a solution of
the hyperbolic Dirac equation (\ref{IBVP0}) and $L:V\to \tilde V$
be a Dirac bundlemorphism. Then $\tilde u=Lu$ is
 a solution of
the hyperbolic Dirac equation (\ref{IBVP0}) with $D$ replaced with $\tilde D$.
Therefore, 
 the Cauchy data sets $C_0(D),\, C_0(\tilde D)$ satisfy
\ba
 K(C_0(D))= C_0(\tilde D),
 \ea
with $K=L|_{W}$.
This shows that the following uniqueness result is optimal.

\begin{theorem}
\label{main1}
Let $D$ and $\tilde D$ be Dirac-type operators  of the 
form (\ref{dirac-type, maximal})
on  the Dirac bundles 
with chirality
$(V,\, F)$ and $(\tilde V,\, \tilde F)$. 
Assume 
%that the induced bundle structures (\ref{b-structure}) on $\p M$ and $\p \tilde M$ 
%are isomorphic, that is, 
that there is an isometry $\kappa:\p M \to
\p \tilde M$ and an  isomorphism 
between the induced bundle structures of $\p M$
and $\p \tilde M$, $K:W\to \tilde W$, which is
compatible with $\kappa$.  Assume also that the Cauchy data sets satisfy
\ba
 K(C_0(D))= C_0(\tilde D).
\ea
Then there is an isometry $\ell: (M,g)\to (\tilde M,\tilde g)$ with
$\ell|_{\p M}=\kappa$ and
a Dirac bundlemorphism $L:V\to \tilde V$ with $L|_W=K$.
In particular, $\tilde D=LDL^{-1}$.
\end{theorem}

%In particular, above theorem concludes that if the boundary data
%are isomorphic then $M$ and $\tilde M$ are isometric and the bundles
%$V$ and $\tilde V$ and Dirac-type operators $D$ and $\tilde D$
%are Dirac bundlemorphic in map $L$ satisfying (\ref{eq: lboro A}).
Note that this theorem states  that the Dirac operators are
equal up to a bundlemorphism, rather than
the connections $\nabla$ and $\tilde \nabla$.

For the inverse boundary spectral problem we obtain a similar
result.

\begin{theorem}
\label{main2}
Let $D_P$ and $\tilde D_{\tilde P}$ be self-adjoint Dirac-type operators  of  
form (\ref{dirac-type})
on  the  Dirac bundles with chirality
$(V,\, F)$ and $(\tilde V,\, \tilde F)$. 
Assume 
%that the induced bundle structures (\ref{b-structure}) on $\p M$ and $\p \tilde M$ 
%are isomorphic, that is, 
that there is an isometry $\kappa:\p M \to
\p \tilde M$, an isomorphism $K:W\to \tilde W$ 
compatible with $\kappa$, and 
that
the eigenvalues and the normalized eigenfunctions of  $D_P$ and $\tilde D_{\tilde P}$ satisfy
\ba
\tilde \lambda_j=\lambda_j,\quad \tilde\phi_j\circ \kappa|_{\p M}= K(\phi_j|_{\p M}),\quad
j\in \Z_+.
\ea
Then there is an isometry $\ell: (M,g)\to (\tilde M,\tilde g)$ with
$\ell|_{\p M}=\kappa$ and
a Dirac bundlemorphism $L:V\to \tilde V$ with $L|_W=K$
such that 
\ba
\tilde D_{\tilde P}= L D_PL^{-1},\quad \tilde P= K PL^{-1}.
\ea
\end{theorem} 

\begin{remark} {\rm To consider above results in more practical terms,
assume that  we  are given the boundary $\p M$, the bundle $W$ on it,
and the Cauchy data set $C_0(D)$ or the boundary spectral data.
The proofs of the above theorems are constructive, 
{\newwtextt providing a method
to 
recover  the manifold $M$ (up to an isometry),
the bundle $V$ with its Hermitian, 
$\bra \cdot,\,\cdot \cet$ and Clifford module, $\gamma$ 
structures, }
and the operators $D$ and $F$
 (up to a Dirac bundlemorphism).}
\end{remark}

%Finally, on the physical interpretation of the results we note that the hyperbolic Dirac-type operator  in (\ref{IBVP0}) is 
%itself time-independent,
%i.e.\ physically the results concern the near-field measurements
%from a stationary scatterer.

\section{Hyperbolic Dirac equation and 
unique continuation}

\noindent {\bf  3.1.\ Basic properties.} We start with 
a general boundary condition $P$
satisfying (\ref{anticommutation}) and (\ref{Fredholm}).
Denote by  $u^f(t)=u^f(x,t)$ the solution of
the initial boundary value problem  (\ref{IBVP0})-(\ref{IBVP0line 2}).

By the  Bochner-Lichnerowicz formula, 
see e.g.\  \cite{Ber,Bs},  
\ba
D^2_0=\nabla^*\nabla+{\cal R},
\ea
where ${\cal R}$ is the curvature endomorphism on  $V$.
Hence $D^2=\nabla^*\nabla+B(x,D)$, where $B(x,D)$ is a first order
differential operator. 
Thus the hyperbolic Dirac-type equation in (\ref{IBVP0})
 implies the wave equation 
\beq\label{wave eq}
(\p_t^2+\nabla^*\nabla+B(x,D))u^f(x,t)=0.
\eeq
%Therefore we sometimes use the name "wave" for solutions $u^f_{\pm}(x,t)$.

In the rest of  subsection 3.1 we consider the
local boundary condition associated to operator 
$P_\Gamma$, cf.\ (\ref{PGamma}).
Let $u=u^f(x,t)$, $t\in \R_+$ be the solution of
\beq\label{IBVP}
(i\p_t+D)u(x,t)&=&0\quad\hbox{in }M\times \R_+,\\
\nonumber
P_{\Gamma} (u|_{\p M\times \R_+})&=&f,\quad
u|_{t=0}=0,
\eeq
 where $f\in \Cnull(\p M \times \R_+,P_{\Gamma} W)$.

Since $D_{\Gamma}$ is a  self-adjoint operator with 
 ${\cal D}(D_{\Gamma}) \subset H^1(M,V)$, the map
$
{\cal U}:f \mapsto u^f
$ 
has an extension 
\ba
{\cal U}\hspace{-1mm}:\hspace{-1mm}H^1_0(\p M \times [0,T],
 P_\Gamma W)\hspace{-1mm}\to
\hspace{-1mm} 
C([0,T],H^1(M,V))\cap C^1([0,T],L^2(M,V)),
\ea 
when $T>0$.

Next we discuss the finite velocity of the wave propagation for the
initial boundary value problem (\ref{IBVP}).
To describe it, let  $\Sigma\subset\partial M$ be open. 
 We define the {\em domain of 
influence} of $\Sigma$ at time $T$ as 
\[ 
 M(\Sigma,T)=\{x\in M\mid \dist(x,\Sigma)\leq T\}, 
\] 
where $\dist(x, y)$ is the distance in  (M, $g)$.  
{\newstuff Using  local trivializations of the Dirac bundle and applying 
classical results on the 
 finite velocity of the wave propagation in hyperbolic systems due to
 e.g. Wilcox \cite{Wilcox}, we see that} 
if $\supp(f) \subset \Sigma \times [0,T]$, then
\beq\label{domains}
\supp (u^f(T)) \subset M(\Sigma,T).
\eeq

Finally, using the local boundary condition, we prove 
Lemma \ref{lem: equivalence of data}.

\noindent {\em Proof:} Assume that we are given $C_0(D)$ and $P$.
Then $\Lambda_P$ is determined by its graph
$\{(Pu,u):\ u\in
C_0(D)\}$
and visa versa. Since the  induced bundle structure
on $\p M$ determines $P_\Gamma$, $C_0(D)$ determines $\Lambda_{\Gamma}$.
\proofbox\smallskip

 \noindent {\bf  3.2.\ Unique continuation.}
Our next goal is to supplement the finite velocity result by
a unique continuation result of Tataru's type,
see \cite{Ta1,Ta3} for the pioneering works. To this end,
for any open $\Sigma \subset \p M$,
we define
  the double cone of influence,
\ba
K_{\Sigma,T}=\{(x,t)\in M\times \R:\ |t-T|+\dist(x,\Sigma)<T\} 
\subset M \times \R.
\ea

\begin{theorem}\label{UCP} 
Let $\Sigma\subset \p M$ be open, $u \in C((0, 2T),\, L^2(M, V))$
be a weak solution of the hyperbolic equation 
$(i\p_t+D)u=0$  in 
$M\times (0,2T)$. Assume that
\beq
\label{3.1.1}
u|_{\Sigma\times (0,2T)}=0.
\eeq
Then $u(x,t)=0$ for $(x,t)\in  K(\Sigma,T).$
\end{theorem}

\noindent {\em Proof:} Assume first that  $u\in C^\infty(M\times (0,2T),\,V)$. 
Condition (\ref{3.1.1}) and the fact that $u$ solves (\ref{IBVP0})
imply that 
\beq
\label{3.1.2}
\p_{N} u|_{\Sigma\times (0,2T)}=0,
\eeq
{\newwtextt where $\p_N$ is the normal derivative. }
By (\ref{wave eq}), $u$ satisfies the wave equation which, in a local 
trivialisation $\Phi:\pi^{-1}(U)\to U\times \C^d$ where
 $U\subset M$ is open, may be written componentwise as
\beq
\label{3.1.3}
(\p_t^2+D^2) u = (\p_t^2 -\Delta_g) I u +{\tilde B}_1(x, \p_x) u =0.
\eeq
Here $\Delta_g$ is the scalar Laplace operator
on $(M, g)$, $I$ is identity matrix,
 and  ${\tilde B}_1$ is 
a first-order operator. We intend to apply Tataru's unique 
continuation result 
\cite{Ta1}, in order to infer from (\ref{3.1.1}), (\ref{3.1.2}), and (\ref{3.1.3})
that $u$ is zero in $K(\Sigma,T)$. 
The difficulty is that the proof in \cite{Ta1} deals only with a scalar wave equation.
In the case of  the scalar wave equation with  time-independent
coefficients a variant of the proof is given in \cite
[Section 2.5]{KKL}. This proof  can be extended, 
word by word, to the vector wave equations of the form
$(\p_t^2 -\Delta_g)I u +{\tilde B}_1 u=0$, where the leading
order term is a scalar operator times the identity matrix. 
Therefore,  (\ref{3.1.1}), (\ref{3.1.2}), and (\ref{3.1.3})
imply that $u=0$ in  $ K(\Sigma,T).$

Let now  $u \in C((0, 2T),\, L^2(M, V))$ be a weak solution 
of the hyperbolic Dirac-type equation  (\ref{IBVP0})
 with the boundary condition 
(\ref{3.1.1}) understood in a weak sense. To prove the claim,
 we start by smoothing 
$u$ with respect to time. Let $\psi(s) \in C^{\infty}_0(-1, 1), \,
\int \psi ds=1$. Then,
\bfo
u_{\e}(x,t) = \int_\R u(x,t-s)\psi_{\e}(s)ds, \quad 
\psi_{\e}(s)= \frac{1}{\e} \psi(\frac{s}{\e}),
\efo
is in $C^{\infty}((\e, 2T-\e), \, L^2(M, V))$ and solves there
the  Dirac-type wave equation 
$(i\p_t+D)u_\e=0$
with  $u_\e|_{\Sigma\times (\e,2T-\e)}=0$. As
\bfo
D^ 2 u_{\e} = - \p_t^2 u_{\e} \in C^{\infty}((\e, 2T-\e), \, L^2(M, V)),
\efo
the standard techniques of the 
elliptic regularity theory, e.g.\  \cite{Gil-Tr}, implies
%\bfo
%D u_{\e} = -i \p_t u_{\e} \in C^{\infty}((\e, 2T-\e), \, L^2(M, V)),
%\efo
%implying 
that $u_{\e} \in  C^{\infty}((\e, 2T-\e), \, H^1(M, V))$.
Iterating this analysis we see that
$
u_{\e} \in  C^{\infty}(M \times (\e, 2T-\e), \, V).
$
By the first part of the proof, 
\ba
%\label{3.1.4}
u_{\e} = 0 \quad  \hbox{for}\,\, \{(x,t): \, \dist(x, \Sigma) + |t-T| < T-\e\}.
\ea
As $u_{\e} \to u$, when $\e \to 0$, this implies 
that $u=0$ in  $ K(\Sigma,T).$
%with a well defined boundary trace, can be handled
%by smoothing $u$ in time. Ind eed,  consider the wave 
%$u_\psi(x,t)=\int_\R u(x,t-s)\psi(s)ds$
%where $\psi\in C^{\infty}_0(-\e,\e), \, \int \psi ds =1$. 
%As $u_\psi$ is also a solution 
%of the hyperbolic Dirac-type equation vanishing on $\Sigma\times (\e,2T-\e)$, 
%we can apply the above unique continuation result to $u_\psi$.
%As $u_\psi \to u$, when 
% $\e \to 0$,
%we prove the unique continuation  for $u$.
\hfill\proofbox

\section{Inner products and controllability.}

\noindent {\bf  4.1.\ Inner products.}
In this section we consider solutions $u^f(t)=u^f(x,t)$
of (\ref{IBVP0})-(\ref{IBVP0line 2}),
 $ f \in P\Cnull(\p M \times \R_+,\, W)$. We compute
the norms of $u^f_+(t)=\Pi_+ u^f(t)$ and $u^f_-(t)=\Pi_- u^f(t)$ using 
the boundary data, that is,
the induced bundle structure
(\ref{b-structure}) on $\p M$ and 
the operator $\Lambda_{P}$.
This type of formulae are called the Blagovestchenskii
formula, due to their first appearance in \cite{Bla} for the one-dimensional 
inverse problems.

Using (\ref{chirality2}) and (\ref{respect})
we see that $F D+DF=0$ and thus
$
\Pi_+ D= D\Pi_-$ and $\Pi_- D= D\Pi_+.
$
%Similarly, since $Q$ satisfies (\ref{respect}), 
%\ba
%\Pi_+ Q= Q\Pi_-,\quad \Pi_- Q= Q\Pi_+.
%\ea
Thus, we can write the hyperbolic Dirac-type equation (\ref{IBVP0}) as
\beq\label{IBVP 1b}
\left(i\p_t
+ \left(\begin{array}{cc} 
0 &D   \\ 
D & 0
\end{array}\right)\right)
 \left(\begin{array}{c}  
u_+(x,t) \\ 
u_-(x,t)
\end{array}\right) =0\quad\hbox{in }M\times \R_+,
\eeq
where $u_+(x,t)=\Pi_+ u(x,t)$ and $u_-(x,t)=\Pi_- u(x,t).$

\begin{theorem}\label{blacho} 
Let $f, h\in P\Cnull(\p M \times \R_+, W)$ be given and 
$u^f, \,u^h$ be the solutions of
 (\ref{IBVP0})--(\ref{IBVP0line 2}).
 Then the response operator $\Lambda_{P}$, and
the induced bundle structure
(\ref{b-structure}) on $\p M$ determine the inner
products  
\beq
\label{25.11.1} 
\bbra  u^f_+(t),\, u^h_+(s)\ccet \quad \hbox{and} \quad
\bbra  u^f_{-}(t), \,u^h_{-}(s) \ccet\quad 
\hbox{for}\ s,t\geq 0.
\eeq
In particular, these data determine
 $\bbra  u^f(t), u^h(s)\ccet$ for $s,t\geq 0$.
\end{theorem}

\noindent {\em Proof:} We consider only the first inner product
in (\ref{25.11.1})
as the second one can be evaluated in the same manner. 
 The Stokes' formula for the Dirac-type operator $D$ takes
the form,
\beq\label{stokes}
\bbra D\phi,\psi\ccet-\bbra \phi,D\psi\ccet=
\int_{\p M}\bra \gamma(N)\phi,\psi\cet dA_g, 
\eeq
where
$
\phi, \psi \in H^1(M,\,V)$ and
 $dA_g$ is the Riemannian volume  on $(\p M,g)$.
Denote  by $ 
I(s,t) = \bbra \Pi_+u^f(t),\Pi_+u^h(s)\ccet.
$
Then, using  (\ref{commutator of gamma}), (\ref{respect}),
 and (\ref{stokes}),
we get
\beq\nonumber 
& & (\partial_s^2 -\partial_t^2)I(s,t) 
 = \bbra \Pi_+ u^f(t),\p^2_su^h(s)\ccet 
- \bbra \Pi_+ \p^2_tu^f(t),u^h(s)\ccet\\
\nonumber
 &=& - \bbra \Pi_+ u^f(t),D(Du^h)(s)\ccet+
\bbra D \Pi_+  u^f(t),Du^h(s)\ccet+
\\
\nonumber
& &-\bbra D\Pi_+  u^f(t),Du^h(s)\ccet+
\bbra  D^2 \Pi_+ u^f(t), u^h(s)\ccet \\
\nonumber
 &=& 
\int_{\p M} \left(\bra \gamma(N)\Pi_+ u^f(t),D \Pi_+u^h(s)\cet
+\bra \gamma(N)D \Pi_+ u^f(t),\Pi_+ u^h(s)\cet \right)\,dA_g
\\
\nonumber
 &=& 
i\int_{\p M} \bigg( \bra \gamma(N)\Pi_+\Lambda_P f(t),\Pi_-\p_s 
(\Lambda_Ph(s))\cet
-\\
\label{A-fomrmu}
 & &\quad \quad \quad \quad \quad \quad - 
\bra \gamma(N)\Pi_- \p_t(\Lambda_P f)(t),\Pi_+ \Lambda_Ph(s)\cet \bigg)\,dA_g.
\eeq
%As $u^f|_{\p M\times \R_+}=P_{\Gamma} u^f + (I-P_{\Gamma}) u^f=
%f-\gamma(N)\Lambda_{\Gamma} f$
%and similar for $u^h$,
%this implies that
%\beq\label{B-fomrmu} 
%& & (\partial_s^2 -\partial_t^2)I(s,t)
%\\ \nonumber
%&=& 
%\frac i4 \int_{\p M}\bra \gamma(N)(1+F)(f-\gamma(N)\Lambda_{\Gamma} f)(t),
%\p_s\{(I-F)(h-\gamma(N)\Lambda_{\Gamma} h)\}(s)\cet\,dA_g 
%\\  \nonumber
%& &
%-
%\frac i4\int_{\p M}\bra \p_t\{(I-F)(f-\gamma(N)\Lambda_{\Gamma} f)\}(t),  
%\gamma(N)(1+F)(h-\gamma(N)\Lambda_{\Gamma} h)(s)\cet\,dA_g.
%\eeq
Since $I(0,t)=I_s(0,t)=I(s,0)=I_t(s,0)=0$,
 we see that $I(s,t)$
satisfies the one-dimensional wave equation (\ref{A-fomrmu})
in the quarter plane $s, t >0$ 
with the known right-hand side
and  homogeneous initial and boundary data. Thus, we can find 
$I(s,t)= \bbra \Pi_+u^f(t),\Pi_+u^h(s)\ccet$ 
for any $s,t \in \R_+$. 
\hfill$\Box$ 

\begin{remark}\label{Blago-remark}{\rm From the proof of Theorem \ref{blacho}, we obtain an
 explicit formula  
 \beq
\label{innerproduct}
\bra \bra \Pi_{\pm} u^f(T),\, u^h(T)\cet \cet ={\mathcal L}_{\pm}[f,h]:=
\int_{\R}\int_{\R} \chi_C(t,s) J_{\pm}^{t,s}(f,h)\, ds dt.\hspace{-2cm}
\eeq
Here, the quadratic form ${\mathcal L}_{\pm}[\cdotp,\cdotp]$ 
is defined using the indicator function  $\chi_C$ is  
of the triangle $C=\{(t,s)\in \R^2:\ T-t\geq |T-s|\}$
and the quadratic form $ J^{t,s}_{\pm}(f,h)$, 
\beq
\label{U-formul.}
\quad \quad J^{t,s}_{\pm}(f,h)&=&
 \frac i2 \int_{\p M} \bigg( \bra \gamma(N)\Pi_\pm \Lambda_P f(t),\Pi_\mp 
\p_s \Lambda_P h(s)\cet-
\\ \nonumber
& &\quad \quad 
-\bra \gamma(N)\Pi_\mp \p_t\Lambda_P f(t),\Pi_\pm \Lambda_Ph(s)\cet \bigg)\,dA_g.
\eeq
}\end{remark}

%\begin{corollary}\label{blacho1} 
%For any $f, h\in P\Cnull(\p M \times \R_+, \, W)$ 
%the operator $\Lambda_{P}$
%and  the
%induced bundle structure
%(\ref{b-structure}) on $\p M$ determine the inner
%product
%\beq
%\label{25.11.1b} 
%\bbra  u^f(t), u^h(s)\ccet \quad \hbox{for}\ s,t\geq 0.
%\eeq
%\end{corollary}
%\HOX{Remove this corollary}
%
%{\em Proof:}  The result clearly follows from Theorem \ref{blacho}. It 
%is, however, useful to give another method to evaluate this
%inner product.
%
%
%Let now $I(s,t) = \bbra  u^f_+(t), u^h_+(s)\ccet$. Then
%\ba
%& &\p_sI +\p_tI = \bbra  u^f(t),\p_su^h(s)\ccet +
%\bbra  \p_t u^f(t),u^h(s)\ccet
%\\
%& &
%=-i \int_{\p M} \bra \gamma(N) \Lambda_P f(t),\Lambda_P h(s)\cet\, dA_g.
%\ea
%As $I(0,t)=I(s, 0)=0, \, t, s \geq 0$, we can find
%$I(s,t)$ from this equation by integrating along the characteristics
%$t-s= \hbox{const}$.
%\hfill$\Box$ 

\noindent {\bf  4.2.\ Global controllability.}
Let  $\Sigma \subset \p M$ be an open set and
\ba
%\label{25.11.3}
& & X^P(\Sigma,T)=
\{u^f(\cdotp,T)\in L^2(M,V):\ f\in PC^{\infty}_0(\Sigma
\times (0, T), \, W)\},
\ea
where $u^f(t)$ is the solution of (\ref{IBVP0})-(\ref{IBVP0line 2})
with the boundary condition associated with the operator $P$.
In the following theorem,
\ba\quad\quad
%\label{radius}
\rad(M, \Sigma) = \max_{x \in M} \dist(x, \Sigma),\quad
\rad(M) =\max_{x \in M} \dist(x, \p M).
\ea

\begin{theorem}\label{global control th} 
1. Let $P$ be a boundary condition satisfying (\ref{anticommutation}) and (\ref{Fredholm}) and $T>2\rad(M)$.  
Then the set $X^{P}(\p M,T)$ is  dense in $L^2(M,V)$.

2. Let $P_\Gamma$ be the local boundary condition
given by (\ref{PGamma}), $\Sigma\subset \p M$ be an open non-empty set,
and $T>2\rad(M,\, \Sigma)$. Then the set $X^\Gamma(\Sigma,T)=X^{P_\Gamma}(\Sigma,T)$ 
is  dense in $L^2(M,V)$.
\end{theorem}

\noindent {\em Proof:} We start by proving  claim {\it 2}.
Let $\eta\in L^2(M,V)$ be orthogonal to $X^{\Gamma}(\Sigma,T)$.
Consider the dual initial boundary value problem 
\beq\label{adjoint IBVP}
(i\p_t+D)v_{\eta}(x,t)&=&0\quad\hbox{in }M\times \R,\\
\nonumber
P_{\Gamma} v_{\eta}|_{\p M\times \R}&=&0, \quad v_{\eta}|_{t=T}=\eta.
\eeq
Integrating by parts we see that,
for any $f \in C^{\infty}_0(\Sigma \times [0,T], \,P_{\Gamma} W )$, 
\beq\label{C-fomrmu} 
0&=&-i\bbra u^f(T),v_{\eta}(T)\ccet
\\
\nonumber
&=& -\int_0^T \int_{ M} \left(\bra i\p_tu^f(t),v_{\eta}(t)\cet-
\bra u^f(t),i\p_tv_{\eta}(t)\cet\right)\,dV_g dt
\\
\nonumber
 &=& \int_0^T \int_{M}\left(\bra Du^f(t),v_{\eta}(t)\cet-
\bra u^f(t),Dv_{\eta}(t)\cet\right)\,dV_g
 dt
\\
\nonumber
 &=& \int_0^T \int_{\p M} \bra \gamma(N)u^f(x,t),v_{\eta}(x,t)\cet
 dA_g(x)dt
\\
\nonumber
 &=& \int_0^T \int_{\p M} 
\left(\bra \gamma(N)P_{\Gamma} u^f(x,t),(I-P_{\Gamma}) 
v_{\eta}(x,t)\cet+\right.
\\ \nonumber
 & &\quad\quad\quad
+
\left.\bra \gamma(N)(I-P_{\Gamma}) u^f(x,t),
P_{\Gamma}  v_{\eta}(x,t)\cet\right)\,
 dA_gdt
%\\
\eeq
\beq
\nonumber
 &=& \int_0^T \int_{\p M} 
\bra \gamma(N)f(x,t),(I-P_{\Gamma}) v_{\eta}(x,t)\cet
 dA_gdt.
\eeq
Here, the map
\beq
\label{28.11.1} 
{\cal B}: L^2(M,V) \mapsto 
{\cal D}'(\p M\times \R,W), \quad {\cal B}(\eta)=
 v_{\eta}|_{\p M\times \R},
\eeq
is continuous.
This may be shown   by an
approximation of $\eta \in L^2(M,\,V)$
by $C^{\infty}_0(M,\, V)$ sections and applying the uniform 
boundedness principle
for the  distributions.  
Since equation (\ref{C-fomrmu}) is valid for any 
$f\in C_0^\infty(\Sigma\times[0,T],P_{\Gamma} W)$,
we see that $(I-P_{\Gamma} ) v_{\eta}|_{\Sigma\times (0,T)}=0.$
Together with the boundary condition of
(\ref{adjoint IBVP}), this makes it possible to apply
Theorem \ref{UCP}. Hence $v_{\eta}$ vanishes in $K_{\Sigma,T/2}$.
However, the set $K_{\Sigma,T/2}$ contains an open neighborhood
of $M\times \{\frac T2\}$. As also 
$P_{\Gamma} v_{\eta}|_{\p M\times \R}=0$, this implies 
that
  $v_{\eta}=0$ in $M\times \R$. In particular,  
 $\eta=v_{\eta}(T)=0$ proving claim {\it 2}.

Claim {\it 1.} \ may be proven by the same arguments as claim {\it 2.}
by taking $P$ instead of $P_\Gamma$ and using property (\ref{anticommutation}).
\hfill\proofbox\smallskip

{
As can be seen in the above proof, for any $t_0$ and 
$f \in C^{\infty}_0(\p M \times (0,t_0), \,P_{\Gamma} W)$,
and $v_\eta$ solving (\ref{adjoint IBVP}) with $\eta\in L^2(M,V)$,
\beq
\label{15.11.1}\quad\quad
\bra\bra u^f(t_0),\,v_{\eta}(t_0)\cet\cet
= i \int_0^{t_0} \int_{\p M} \bra \gamma(N)f(x,t),\, v_{\eta}(x,t) \cet \,dA_g dt,
\eeq
where the right hand side is understood as a distribution pairing.

%Observe also that $\Lambda_+,\, F|_{\p M},\,\gamma(N)|_{\p M}$ 
%determine the Cauchy data
%\ba
% f \mapsto u^f|_{\p M \times \R_+}= (f,\, -\gamma(N) \Lambda_+f),
%\quad f \in C^{\infty}_0(\p M \times \R_+,\, W_+),
%\ea
%and, therefore, determine the operator $\Lambda_-$.

\smallskip

\noindent {\bf  4.3.\ Local controllability results.}
When $P=P_\Gamma$ is the local projection on $\p M$, we define
\ba
%\label{25.11.3B}
& &X^\Gamma_+ (\Sigma,T)=\{\Pi_+ u^f(\cdotp,T)
\in L^2(M,V_+):\ f\in C^{\infty}_0(\Sigma
\times (0, T), \,P_{\Gamma} W\},\\
\nonumber
& &X^\Gamma_- (\Sigma,T)=\{\Pi_- u^f(\cdotp,T)
\in L^2(M,V_-):\ f\in C^{\infty}_0(\Sigma
\times (0, T), \, P_{\Gamma} W)\}.
\ea
For an arbitrary set $S\subset M$,  
we also denote 
$$L^2(S,\,V)=\{v\in L^2(M,V):\supp(v)\subset 
\overline S\}.$$

\begin{theorem}\label{local control th} 
Let $P_\Gamma$ be the local projection on $\p M$
given by (\ref{PGamma}). Then, for any   $T>0$ and open $\Sigma\subset \p M$, 
the sets $X^\Gamma_{ \pm} (\Sigma,T)$ are dense in the subspaces
$\Pi_{\pm} L^2(M(\Sigma,T),\,V) \subset L^2(M, V)$, correspondingly. 
\end{theorem}

\noindent {\em Proof:}  
{We will prove this result for $X^\Gamma_+ (\Sigma,T)$, the
claim for  $X^\Gamma_- (\Sigma,T)$ may be  obtained in the same manner.

Let $\eta \in \Pi_+L^2(M(\Sigma,T),V)$ 
satisfy
\bfo
\bra\bra \eta,\, \Pi_+u^f(T) \cet\cet =0 \quad \hbox{for any} \,\, 
f \in C^{\infty}_0(\Sigma \times (0, T), \,P_{\Gamma} W).
\efo
Consider the initial boundary value problem
\beq\label{adjoint IBVP2}
(i\p_t+D)v_{\eta}(x,t)&=&0\quad\hbox{in }M\times \R, \\
\nonumber
P_{\Gamma} v_{\eta}|_{\p M\times \R}&=&0, \\
\nonumber
\Pi_+ v_{\eta}|_{t=T}=\eta,& & \Pi_- v_{\eta}|_{t=T}=0.
\eeq
Then $v_{\eta} \in C(\R,\,L^2(M,V))$ and, {\newwtext 
as noted in the proof of Theorem \ref{global control th}, it has
a well defined trace}
$v_{\eta}|_{\p M \times \R} \in {\mathcal D}'(\p M \times \R)$. Same considerations 
as in (\ref{C-fomrmu}) show that
\bfo
0 = \int_0^T \int_{\p M} 
\bra \gamma(N)f(x,t),\, (I-P_{\Gamma}) v_{\eta}(x,t)\cet\,
dA_gdt,
\efo
for any  $f\in C^\infty_0(\Sigma\times[0,T],P_{\Gamma} W)$. Therefore, 
$(I-P_{\Gamma}) v_{\eta}|_{\Sigma \times (0,T)}=0$.
Together with the boundary condition in 
(\ref{adjoint IBVP2}), this imply that $v_{\eta}|_{\Sigma \times (0,T)}=0$, i.e.
$\Pi_{\pm}v_{\eta}|_{\Sigma \times (0,T)}=0$.  Using  (\ref{IBVP 1b}) 
and the fact that $\Pi_-v_{\eta}|_{t=T}=0$,
we see that $\Pi_-v_{\eta}(T+s)= -\Pi_-v_{\eta}(T-s)\,$
and $ \Pi_+v_{\eta}(T+s)= 
\Pi_+v_{\eta}(T-s)$. 
Therefore, 
\beq
\label{63a}
\Pi_{\pm}v_{\eta}|_{\Sigma \times (T,2T)}=0 \quad \hbox{ and} \quad 
\hbox{supp}\left( v_{\eta}|_{\Sigma \times (0,2T)} \right) 
\subset \Sigma \times\{T\}.
\eeq
{We intend to show that $v_{\eta}|_{\p M \times (0,2T)}=0$. Due to the boundary condition in
(\ref{adjoint IBVP2}) this is equivalent to
\beq
\label{additional1}
\int_0^{2T} \int _{\Sigma}\bra\gamma(N) f(x,t),\, 
v_{\eta}(x,t)\cet \,dA_g dt=0
\eeq
for any $f\in C^\infty_0(\Sigma\times(0,2T),P_{\Gamma} W)$. In view of (\ref{63a}),  it is 
enough to prove that
 \beq
 \label{15.11.2}
 & & \int_{T-2\e}^{T+2\e} \int _{\Sigma} 
\bra \gamma(N) f_{\e}(x,t),\,v_{\eta}(x,t)\cet \,dA_g dt 
\\ \nonumber
& &=
i\bra\bra u^{f_\e}(T+2\e),\,v_{\eta}(T+2\e)\cet \cet \longrightarrow 0,
\quad \hbox{when} \,\,\e \to 0.
 \eeq
Here $f_\e(\cdot,t) = \chi_\e(t-T) f(\cdot,t)$ and $\chi_\e(s)$ 
is 
a smooth cut-off function,
 $\chi(s)=1$ for $s \in (-1,1)$,
  $\,\hbox{supp}(\chi) \subset (-2,2)$ with
  $\chi_\e(s) = \chi(s/\e)$.
  Note that the equality in (\ref{15.11.2}) follows from (\ref{15.11.1}).
  }
%  Thus, to prove that $v_{\eta}|_{\Sigma \times (0,2T)}=0$ it is 
%sufficient to show that the right hand side in (\ref{15.11.2}) 
%tends to $0$ when $\e \to 0$ for any 
%$f\in C^\infty_0(\Sigma\times(0,2T),W_+)$.

 Assume, in the beginning, that 
$f_\e(x,t) = \chi_{\e}(t-T) F_{\p}(x)$, where
 $ F_{\p}= F|_{\p M}, F \in C^{\infty}(M,V)$. We 
can represent the solution as
 \bfo
 u^{f_\e}(x,t)= \chi_{\e}(t-T) F(x) + \omega_\e(x,t),
 \efo
 where $\omega_\e$ is the solution to
 \ba
& &\hspace{-5mm} (i \p_t +D) \omega_\e = -(i \p_t +D)\,(\chi_\e F )= 
-\frac{i}{\e} \chi'\left((t-T)/\e\right) F -
  \chi_\e(t-T) \,D F,
 \\
& &\hspace{-5mm}  P_{\Gamma} \omega_\e|_{\p M \times \R}=0,
 \quad
 \omega_\e|_{t=T-2\e}=0.
 \ea
We represent  $F$ and $D F$ as $L^2(M)$-converging series of
eigenfunctions $\phi_k$ of $D_{\Gamma}$,
 $$
 F= \sum_{k=1}^\infty a_k \phi_k,\quad  D F = \sum_{k=1}^\infty b_k\phi_k.
 $$
 Then
 \beq
 \label{15.11.3}\quad\quad
 \omega_\e(T+t)= \sum_{k=1}^\infty \left(w^1_k(t)+w^2_k(t)\right)\phi_k = 
\omega_\e^1(T+t)+\omega_\e^2(T+t), 
 \eeq
 \[
 w^1_k(t)= \frac{a_k}{i\e} \int_{T-2\e}^t 
e^{i \lambda_k (t-s)} \chi'(s/\e)ds, \ \  
 w^2_k(t)=- b_k\int_{T-2\e}^t e^{i \lambda_k (t-s)} \chi (s/\e)ds.
\]
  It then follows from
 (\ref{15.11.3}) that
  \beq
 \label{15.11.4}\quad\quad
 \|\omega_\e^1(T+2\e)\|_{L^2} \leq C_1 \|F\|_{L^2}, 
\quad \|\omega_\e^2(T+2\e)\|_{L^2} \leq C_2 \e  \|D F\|_{L^2}.
 \eeq
 What is more, introducing
 \bfo
 \omega^{1,N}_\e(T+t) = \sum_{|\lambda_k| >N} w^1_k(t) \phi_k,
 \efo
 we see that 
  \beq
 \label{15.11.5}\quad \quad
 \lim_{N\to \infty}\|\omega^{1,N}_\e(T+2\e)\|_{L^2}=0,
 \eeq 
 uniformly with respect
 to $\e \in (0,\e_0)$. To estimate $\omega_\e^1-\omega_\e^{1,N}$, 
 use integration by parts in the equation for $w^1_k(t)$ in
 (\ref{15.11.3}) to get
  \beq
  \label{15.11.6}
 |w^1_k(2\e)| \leq C_3N |a_k| \e\quad\hbox{when }|\lambda_k|\leq N.
 \eeq
 Combining (\ref{15.11.5})--(\ref{15.11.6}) we see that
$\|\omega_\e^{1}(T+2\e)\|_{L^2}\to 0$ when 
$\e\to 0$. This and  the estimate for 
$\omega_\e^{2}(T+2\e)$ in (\ref{15.11.4}) imply the desired
 result (\ref{15.11.2}) for $f_\e(x,t) = \chi_{\e}(t-T) F_{\p}(x)$.

 Returning to a general $f \in C^{\infty}_0(\Sigma\times (0,2T),\, P_{\Gamma} W)$,
 it remains to estimate $\omega_\e(T+2\e)$
 for $f_\e(x,t) = (t-T)\chi_{\e}(t-T) \tilde{f}(x,t)$ 
with smooth $\tilde{f}(x,t)$.
 Arguments similar to those leading to the estimate of 
$\,\omega_\e^{2}(T+2\e)\,$ in (\ref{15.11.4})
 show that $\,\|\omega_\e(T+2\e)\|_{L^2} \to 0$ when $\e \to 0$.
} This proves (\ref{additional1}) and, therefore, Theorem
\ref{local control th}.
\hfill\proofbox\smallskip

We complete this section by describing a procedure to find $\rad(M, \Sigma)$.
\begin{corollary}
\label{rad} For any open $\Sigma \subset \p M$, the Cauchy data set
$C_0(D)$ and the induced bundle structure on $\p M$ determine 
$\rad(M,  \Sigma)$.
\end{corollary}
\noindent {\em Proof:}
By the
%using
%Unique Continuation Theorem \ref{UCP}
 local controllability result, Theorem
\ref{local control th}
we see that, if $T_1>\rad(M, \Sigma)$, then the following
claim $(A)$ is valid: 
\smallskip

\noindent
$(A)$ For any $T_2>0$ and $f\in C^\infty_0(\p M \times \R_+,P_\Gamma W)$
there are $ h_j^+,h_j^-\in C^\infty_0(\Sigma \times (0,T_1),P_\Gamma W)$, $j=1,2,\dots$
such that
\beq
\label{codn A}
& &\lim_{j\to \infty} \| \Pi_+(u^{h_j^+}(T_1)-u^{f}(T_2))\|_{L^2}=0,\\
& & \nonumber
\lim_{j\to \infty} \| \Pi_-(u^{h_j^-}(T_1)-u^{f}(T_2))\|_{L^2}=0.
\eeq

On the other hand, if  $T_1<\rad(M, \Sigma)$, we see that for $T_2> \rad(M, \Sigma)$ there is
$f$ such that $\supp(u^{f}(T_2))\not\subset M(\Sigma,T_1)$. Therefore,
claim (A) is not valid. 

As by Theorem \ref{blacho}
 the norms on the left-hand sides of (\ref{codn A}) can be evaluated
using the boundary data,
we find $\rad(M, \Sigma)$ by taking the  infimum of all $T_1>0$
such that $(A)$ is valid.
\hfill\proofbox

\section{Generalized boundary sources}
\label{sec: Gen. sources}

\noindent {\bf  5.1.\ The wave norm.} We define the wave operators 
\ba
%\label{30.11.3a}\quad\quad
{\mathcal  U}^T:PC^\infty_0(\p M \times (0,T), W )
\to L^2(M,\,V), \quad
{\mathcal  U}^Tf =u^f(T),
\ea 
where $T=3\rad(M)$ and 
$u^f(t)=u^f(x,t)$ is the solution 
of (\ref{IBVP0})-(\ref{IBVP0line 2}). By means of the wave operators, 
we 
define a semi-norm, or the wave semi-norm on the space of 
boundary sources as 
\begin{equation}\label{f norm} 
 \|f\|_{\F} = \|{\mathcal  U}^T f\|_{L^2(M,V)}. 
\end{equation} 
By 
Theorem \ref{blacho}, the knowledge 
of the operator ${\Lambda}_P$ and 
the induced bundle structure
(\ref{b-structure}) on $\p M$
enables 
us to calculate explicitly the semi-norm  (\ref{f norm})
of 
any $f \in P C^\infty_0(\p M \times (0,T), W )$. 
 
 To complete the space of boundary sources, 
we define 
the equivalence relation $\sim$ on this space by setting 
\[ 
 f\sim h\mbox{ if and only if } u^f(T)=u^h(T). 
\] 
Further, we define the space ${\mathcal B}$ as 
\ba %\label{space cal F}
 {\mathcal B}=
P C^\infty_0(\p M \times (0,T), W )/\sim. 
\ea
Then  (\ref{f norm}) becomes a norm on ${\mathcal B}$.
Finally, we complete ${\mathcal B}$  
with respect 
to the norm (\ref{f norm}). This space, 
denoted by $\F= {\overline{\mathcal B}}$, 
consists of the
 sequences of sources which are the Cauchy sequences 
with respect to the norm 
(\ref{f norm}). We denote these sequences by ${\hat f}$,
\[ 
\hat f = (f_{j})_{j=0}^\infty,\quad f_j\in 
PC^\infty_0(\p M\times (0, T), W ). 
\] 
Also, we denote $J:PC^\infty_0(\p M\times (0, T), W )\to \F$ the factor-projection
that maps $f$ to its equivalence class in ${\mathcal B}$.

The sources $\hat f\in \F$ are called {\em generalized sources} 
in the sequel. The corresponding waves are denoted 
by
\beq
\label{gener}
  {\mathcal  W}\hf(t)  = 
\lim_{j\to \infty} u^{f_j}(t)\quad
\hbox{for }t\geq T,
\eeq
where the convergence in the right-hand side is the 
$L^2(M,\,V)-$conv\-er\-gence.
By  the global controllability, Theorem
\ref{global control th}
and the definition (\ref{f norm}) of the norm, the map
$\hat f\mapsto   {\mathcal  W}\hf(t)$ is an isometry,
$\F\to L^2(M,V)$.
 We extend the notation ${\mathcal  W}{\hat f}(t)$ 
for all $t\in \R$ by defining
\beq
\label{semigroup}
 {\mathcal  W}{\hat f}(t)=\exp(i(t-T)D_P){\mathcal  W}{\hat f}(T)
\quad
\hbox{for }t\in \R.
\eeq
Due to this formula we call ${\mathcal  W}\hat f(t)$ 
a {\it generalized
wave.}

\begin{remark} {\rm
The choice of the parameter $T=3\rad(M)$ in the above construction
can be replaced with any fixed $T >\, 2\rad(M)$.
}\end{remark} 

%\begin{remark} {\rm
%It is clear from (\ref{space cal F}) that ${\mathcal B}={\mathcal B}(T)$
%and its completion ${\mathcal F}={\mathcal F}(T)$  depend on $T$.
%However, it follows from
% the global controllability theorem
%\ref{global control th} that all spaces ${\mathcal F}(T)$, $T>2\rad(M)$
%are unitary equivalent to $L^2(M,V)$. Therefore, all  ${\mathcal F}(T)$, $T>2\rad(M)$
%are unitary equivalent and we can identify elements
%$\hat f\in {\mathcal F}(T)$ and $\hat f\in {\mathcal F}(T')$ if
%${\mathcal W}\hat f(t)={\mathcal W}\hat f'(t)$ for $t\in \R$. Thus we 
%can speak about the Hilbert space of the generalized sources $\F$
%without referring to the index $T$. 
%}
%\end{remark}

We note that the above construction of the space of  generalized
sources is well-known in PDE-control, e.g.~\cite{LTr,Ru}. 

For $t \in \R$ and any boundary source 
${\hat f} \in {\F}$,
we can find if the condition
\beq
\label{gener2}
 \p_t^l {\mathcal W}{\hat f}(t) 
%= 
%\lim_{j\to \infty} {\mathcal W}{\p_t^l f_j}(t) 
\in L^2(M,\,V) 
\eeq
is satisfied for a given $ l \in \Z_+$.
Indeed, for $l=1$, condition (\ref{gener2}) 
for $\hat f=(f_j)_{j=1}^\infty\in {\F}$
is equivalent to the 
existence of $\hat h \in {\F}$ 
such that
\beq
\label{F1}
& &
\lim_{\delta\to 0} \|\frac{{\mathcal  W}f(t)-{\mathcal  W}f(t-\delta)}{\delta}-{\mathcal  W}h(t)\|_{L^2(M,V)}
\\
\nonumber
& &
=
\lim_{\delta\to 0} \|
 {\mathcal  W}\hspace{-1mm}\left(\frac {\hat f-{\cal T}_\delta\hat f}\delta - \hat  h\right)
 \hspace{-1mm}(t)\|_{L^2(M,V)}
 =0,
\eeq 
$t >T$.
Here ${\cal T}_\delta: f(x,t) \mapsto f(x,t-\delta)$
is the time-shift
{\newwtextt naturally extended from $PC^\infty_0(\p M\times (0, T), W )$ to
${\mathcal F}$. } 
Note that the validity of (\ref{F1}) is independent of value of $t>T$ as the norms
$\|{\mathcal W}{\hat f}(t)\|_{L^2}$ 
 are independent of  $t$.
Using the induced bundle structure
on $\p M$ and the response operator $\Lambda_{P}$
we can verify 
if a given generalized source $\hat f= (f_j)_{j=1}^\infty$ 
{\newwtextt and a trial $\hat h = (h_j)_{j=1}^\infty$ satisfy}
 condition (\ref{F1}).
In this case we denote
\bfo
{\hat h}= \p_t {\hat f}\quad\hbox{and}\quad \hat f\in {\mathcal D}(\p_t).
\efo
%Note that as $\p_t u^f(t)=u^{\p_t f}(t)$
%in this case $\hat h=(\p_t f_j)_{j=1}^\infty$ gives \HOX{check}
%an explicit representation for $\hat h\in \F$.
Having defined $\p_t^{k-1} {\hat f}$ we  define 
$\p_t^{k} {\hat f}$, if it exists, by induction. The spaces 
of  generalized sources which have $k \in \Z_+$ time derivatives
will be denoted  by $\F^k={\cal D}(\p_t^k)$ with
$\F^\infty=\bigcap_{k\in \Z_+} \F^k$.
Clearly,
\ba
\p_t Jf=J(f_t),\quad f\in PC^\infty_0(\p M\times (0, T), W ), 
\ea
that is, $\p_t$ is an extension of the time derivative
$f\mapsto f_t$ defined in the usual way for $f\in  PC^\infty_0(\p M\times (0, T), W ).$
Note that by (\ref{semigroup}) and (\ref{gener2}),
 \beq
\label{equivalence1}
{\hat f} \in {\cal D}(\p_t^k), \quad \hbox{if and only if} \quad
{\mathcal W}{\hat f}(t) \in {\cal D}(D_{P}^k), \quad t \in \R.
\eeq
In addition, 
 \beq
\label{equivalence2}
{\mathcal  W}(\p_t^k \hat f)(t)= i^k D_P^k {\mathcal W}{\hat f}(t)\quad
\hbox{ for }{\hat f} \in \F^k.
\eeq

\smallskip

\noindent {\bf  5.2.\ Boundary sources corresponding to eigenfunctions.}
When  $\hf \in {\mathcal  F}^1$
 we can evaluate, using Theorem \ref{blacho}, 
the  inner products
\ba
%\label{1-form}
& &\bbra D_P {\mathcal W}{\hat f}(T),\, {\mathcal W}{\hat f}(T)\ccet=
-i\bbra {\mathcal W}(\p_t\hat f)(T)
,\, {\mathcal W}{\hat f}(T)\ccet,\\
& &\bbra D_P {\mathcal W}{\hat f}(T),\, D_P {\mathcal W}{\hat f}(T)\ccet=
\bbra {\mathcal W}(\p_t\hat f)(T)
,\,{\mathcal W}(\p_t\hat f)(T)\ccet,
\ea
as well as the inner products 
$\bbra {\mathcal W}{\hat f}(T),\,{\mathcal W}{\hat f}(T)\ccet$.
On the other hand, 
by the Hilbert-Courant min-max principle 
the eigenvalues $[\lambda_j]^2$ of $D_{P}^2$ can be found as
\beq\label{r-q}
[\lambda_j]^2=
\sup_{u_1,\dots,u_{n-1}}\inf_{u_n\perp u_1,\dots,u_{j-1}}
\frac {\bbra D_P u_j,\, D_P u_j\ccet}{\bbra u_j,u_j\ccet},
\eeq
where $u_j$ are in ${\cal D}(D_{P})$. 
{\newwtextt This makes it possible to prove the following result: }

\begin{proposition}
\label{5.1}
The induced bundle structure
(\ref{b-structure}) on $\p M$
and the response operator $\La_P$
 determine 
the eigenvalues $\lambda_j$ of $D_{P}$  
and the generalized sources 
$\hh_j \in \F^\infty,\,$ $j=1,2,\dots$, 
such that  
$\phi_j={\mathcal W}{{\hat h}_j}(T)$ are orthonormal 
eigenfunctions of $D_{P}$ that form a complete basis of $L^2(M,V)$. 
\end{proposition}

\noindent {\em Proof:} 
Using (\ref{r-q}), we can find all the eigenvalues 
of $D_{P}^2$. Indeed, by the global controllability, Theorem
\ref{global control th} and relation (\ref{equivalence2}), ${\mathcal W}
\hat f(T)$ runs over 
${\mathcal D}(D_P)$ when
$\hat f$ runs over $\F^1$.
At the same time, the minimizers
in (\ref{r-q}) are those generalized
sources  $\hf_{j} \in {\mathcal  F}^1$ that the  waves 
$u_j={\mathcal W}{\hf_j}(T)$ are the corresponding orthonormal
eigenvectors of $D_{P}^2$.
Furthermore, since we can evaluate $\bbra D_{P} u_j,\, u_j\ccet$, 
we can decompose the eigenspace of $D_{P}^2$, which
 corresponds to an eigenvalue $[\lambda_j]^2$, into 
a direct sum of the
eigenspaces of $D_{P}$ which corresponds to the 
eigenvalues $\lambda_j$ and 
$-\lambda_j$.
 Hence, we can find 
$\hh_j \in {\mathcal  F}^1$ such that
$\phi_j={\mathcal W}{{\hat h}_j}(T)$ are the orthonormal 
eigenfunctions of $D_{P}$. Clearly, $\hh_j \in {\mathcal  F}^\infty.$
\hfill\proofbox \medskip
 
%{\newwtext In the following we assume that the eigenvalues are 
%enumerated
%so that $|\lambda_j|\leq |\lambda_{j+1}|$, $j=1,2,\dots$.}

Applying Theorem \ref{blacho} to ${\mathcal W}{{\hat f}}(t)$ 
and ${\mathcal W}{{\hat h}_j}(T)$, we can evaluate the Fourier 
coefficients ${w}_j(t)=\bra\bra {\mathcal W}{{\hat f}}(t),\phi_j\cet\cet$,
$j=1,2,\dots$, of ${\mathcal W}{{\hat f}}(t)$. Using $w_j(t)$,
we can compute the norms
 $\|(1+|D_{P}|)^s {\mathcal W}{\hat f}(t)\|_{L^2(M)}$,
 $s\geq 0$. 
We define the spaces $\F^s$ as  the spaces of those  
generalized sources  $\hat f\in \F$ 
for which
${\mathcal W}{\hat f}(t)\in {\cal D}(|D_{P}|^s)$ and define
\bfo
\|\hat f\|_{\F^s}=
\|\,(1+|D_{P}|)^s {\mathcal W}{\hat f}(t)\|_{L^2(M)}
= \sum_{j=1}^\infty (1+|\la_j|)^{2s} |w_j(t)|^2,
\efo
where the right-hand side is independent  of $t\in \R$.

\section{ Index formulae using boundary data}
In this section we use boundary data to verify if the
 decomposition condition (\ref{decomp 1.}) is valid
and give formulae for the Fredholm index of the operator $ D_P^+$
appearing in decomposition
(\ref{decomposition0}).

\smallskip

\noindent {\bf  6.1.\ Decomposition of $D_P$.}
Let us consider first the Dirac operator $D_\Gamma$ associated
with $P_\Gamma$. Using  properties (\ref{chirality2}) of the 
chirality operator, we can see
that, if 
$u \in  {\mathcal D}(D_{\Gamma})$ and $\Pi_+u \in  {\mathcal D}(D_{\Gamma})
$, then $u|_{\p M}=0$.
Therefore,  the decomposition condition (\ref{decomp 1.}) is
not satisfied for the Dirac-type operator $D_\Gamma$.
However, there is a large number of
important examples of projectors  $P$ 
when  the decomposition condition (\ref{decomp 1.})
is valid for the corresponding Dirac-type operator $D_P$. Some of these are discussed below.

\smallskip
 
\noindent {\bf  6.2.\ Decomposition of $\F$.}
Assume that we are given the operator $\Lambda_{P}$
with a possibly non-local boundary operator $P$ and
the induced bundle structure
 on $\p M$.   
Then the formula  
 (\ref{innerproduct})
determines the non-negative sesquilinear forms 
${\mathcal L}_+[\cdotp,\cdotp]$ and ${\mathcal L}_-[\cdotp,\cdotp]$  
on $PC^{\infty}_0(\p M \times (0, T), W)$ such that
\ba
%\label{B} 
{\mathcal L}_{\pm} (f,h)= \bra \bra \Pi_{\pm} u^f(T),\, u^h(T)\cet \cet. 
\ea
Observe that the quadratic form 
${\mathcal L}={\mathcal L}_++{\mathcal L}_-$
defines a semi-norm $\sqrt{{\mathcal L}[f,f]}$  that coincides,
on 
$PC^{\infty}_0(\p M \times (0, T), W)$,
with the semi-norm defined by (\ref{f norm})
\beq
\label{B slava} 
{\mathcal L} [f,f]={\mathcal L}_+ [f,f]+{\mathcal L}_- [f,f]=\|f\|_\F^2.
\eeq
Therefore, ${\mathcal L}_+$ and ${\mathcal L}_-$ can be extended to $\F$.

In turn, these forms give rise to two bounded  non-negative 
operators 
${B}_{\pm} $ on $\F$ so that
\ba
{\mathcal L}_\pm[\hat f,\hat h]=({B}_{\pm} \hat f,\hat h)_\F.
\ea
The fact that $\Pi_+ L^2(M,V)\oplus \Pi_- L^2(M,V)=L^2(M,V)$ implies
that ${B}_{\pm} $ are orthoprojections in $\F$ satisfying
${B}_+ \oplus  {B}_-=I$. Introducing an operator $S$,
\bfo
S\hat f={\mathcal W }\hat f(0),\quad
\hat f\in \F,
\efo
we have
\beq\label{S-identity}
S{B}_{\pm} \hat f=\Pi_\pm Sf.
\eeq

\smallskip

\noindent {\bf  6.3.\ Index formulae.} At first we  prove Theorem \ref{main0}.

\noindent {\em Proof:} 
{\newwtextt Using the operator $S$, equation (\ref{equivalence2})
may be written as}
%As for any $t\in \R$, $\exp(itD_P)u\in 
% {\mathcal D}(D_P)$ if and only if $u\in 
% {\mathcal D}(D_P)$, it follows from 
% (\ref{semigroup}),
% (\ref{equivalence2}), and (\ref{Slava 87}) that
%$S:\F^k\to {\mathcal D}(D_P^k)$ is an isomorphism and
 \beq\label{Slava 87}
 D_P=S(-i\p_t)S^{-1}.
 \eeq
Together with  (\ref{S-identity}), this equation imply that
 the condition $\Pi_\pm {\mathcal D}(D_P)\subset {\mathcal D}(D_P)$
 is satisfied if and only if ${B}_\pm {\mathcal D}(\p_t)\subset {\mathcal D}(\p_t).$
 This proves claim {\it 1.}
 
If decomposition (\ref{decomp 1.}) is valid, so that $D_P$ can be decomposed
as in (\ref{decomposition0}), formula (\ref{Slava 87}) implies
that $\p_t$ can be decomposed as
\ba
\p_t = \left(\begin{array}{cc} 
0 &\p_t^-   \\ 
\p_t^+ & 0
\end{array}\right): {B}_+{\mathcal D}(\p_t)
\oplus 
  {B}_-{\mathcal D}(\p_t)
\to 
  {B}_-\F
\oplus
 {B}_+\F,
\ea
where  $\p_t^\pm$ are the restrictions of $\p_t$ on $ {B}_{\pm}{\mathcal D}(\p_t)$,
correspondingly.
As $D_P^\pm$ are Fredholm operators,  we conclude from (\ref{S-identity}) and (\ref{Slava 87})
that $\p_t^\pm$
  are Fredholm operators and $\ind\,(\p_t^+)=\ind\, (D_P^+)$. This proves claim {\it 2.}
 \proofbox\smallskip

\begin{corollary}
\label{cor:index}
Assume that condition (\ref{decomp 1.}) is valid for $D_P$.
Then
\bfo
\ind\,(D_P^+)= \hbox{sgn} \left( [B_+ -
B_-]|_ {{\rm Ker}\,(\p_t)}\right).
\efo
\end{corollary}

\noindent {\em Proof:} As $FD+D F=0$, 
{\newwtextt we have $F:\, \hbox{Ker}(D_P)\to \hbox{Ker}(D_P)$ so that }
$\Pi_\pm:\hbox{Ker}(D_P)\to \hbox{Ker}(D_P)$. Thus, $ \ind\,(\p_t^+)$ is equal to
 \bfo
\hbox{dim}\,\hbox{Ker}\,(\p_t|_ {B_+\F^1} )-
\hbox{dim}\,\hbox{Ker}\,(\p_t|_ {B_-\F^1} )=
\hbox{sgn} \left( [B_+ -
B_-]|_ {{\rm Ker}\,(\p_t)}\right),
\efo
{\newwtextt where we use the fact that $B_{\pm}$ are orthoprojectors in ${\mathcal F}$. }
 \proofbox\smallskip

By Proposition
\ref{5.1}, using the boundary data it is possible to find generalized sources ${\hat h}_j
=(h_j^k)_{k=1}^\infty
\in \F^{\infty}$ such that
$
{\mathcal W}{\hat h}_j(T)= \phi_j,
$
with $\{\phi_j\}_{j=1}^{\infty}$ being an orthonormal basis
of the eigenfunctions of $D_P$. Let $\phi_j, \, j=1, \dots, \nu,$
where
$\nu= \hbox{dim}\,\hbox{Ker}\,(D_P)$,
be the eigenfunctions spanning $\hbox{Ker}\,(D_P)$.
Then, for any $t_0>T$,
\bfo
\lim_{k\to \infty} \Lambda_P (h_j^k)|_{\p M\times (t_0,\infty)}
 = 1(t)\, \phi_j(x)|_{\p M}, \quad 
 j=1, \dots, \nu,
\efo
with $1(t)\equiv 1$.
Thus, $\Lambda_P$ determines the subspace
\ba
%\label{0-space}
{\mathcal N}= \hbox{span} \{ \phi_j|_{\p M}: \ 
 j=1, \dots, \nu\}\subset (1-P)C^\infty(\p M,W).
\ea
As
$F:{\mathcal N}\to {\mathcal N}$ and
 $F=\Pi_+-\Pi_-$,   Corollary \ref{cor:index} implies the following result:

\begin{lemma}
\label{index}
Assume that $D_P$ satisfies condition (\ref{decomp 1.}). Then
\bfo
\ind (D_P^+) =  \hbox{sgn} \left(F|_{\mathcal N} \right),
\efo
where the right-hand side can be evaluated in terms of 
the induced bundle structure
(\ref{b-structure}) on $\p M$ and
the response operator $\Lambda_P$.
\end{lemma}

  \noindent {\bf  Example 1, continued.} It follows from the definitions of
  ${\mathcal F}$ and the relative and absolute boundary conditions,
that the Dirac operators $(d+\delta)_r$ and $(d+\delta)_a$
can be decomposed into $(d+\delta)_{r}^{\pm}$
 and $(d+\delta)_{a}^{\pm}$. 
Clearly,
$\hbox{Ker}[(d+\delta)_{r}^{+}]$
 consists of the 
relative harmonic forms of the even order
and $\hbox{Ker}[(d+\delta)_{r}^{-}]$  
 of relative harmonic forms of the odd order. Similarly,  
$\hbox{Ker}[(d+\delta)_{a}^{+}]$
and $\hbox{Ker}[(d+\delta)_{a}^{-}]$ consists of 
absolute harmonic forms of an even and odd order, correspondingly.
Therefore, $\ind[(d+\delta)_{r}^{+}]$ is the Euler characteristic
$\chi(M)$ of $M$. By the Poincare duality,
$\ind[(d+\delta)_{a}^{+}]=\pm \ind[(d+\delta)_{r}^{+}]$
depending on $n=\dim(M)$.
% and
%$\ind (d+\delta)_a=\chi(M)-\chi(\p M)$. \HOX{Check} 
%respectively.
Thus, Theorem \ref{main0}
 gives a representation of the Euler
characteristic $\chi(M)$ in terms of the Cauchy data set $C_0(D)$ of the 
hyperbolic form-Dirac equation.

\smallskip

  \noindent {\bf  Example 2, continued.} Similarly to example 1, the index
of the extended Maxwell operator $i(d-\delta_\alpha)$ 
with the electric boundary condition
is the Euler characteristic $\chi(M)$.
\smallskip

  \noindent {\bf  Example 3, continued.} 
The generalized  
signature operator on even dimensional
manifold with boundary satisfies
decomposition condition (\ref{decomp 1.}). This was
 first considered in the famous paper 
\cite{APS} (for further developments in this area see e.g.\ 
 \cite{Br, gilkey1, gilkey2,  Me1,Muller}). Theorem \ref{main0}
provides a possibility to find the index of the 
 signature operator from the Cauchy data of its Dirac operator. 
%For the the geometrical meaning of
%the index, see e.g.\ \cite{gilkey1}.\HOX{Check reference}

\section{Reconstruction of the manifold}

In this section we will show how to
 reconstruct the manifold, $M$ and metric, $g$ from the
response operator $\Lambda_{\Gamma}$.
To this end, we will show that the boundary data determine the 
set of the {\em boundary distance functions} that 
determine $(M,g)$ up to an isometry.
%The results of this and the next section generalize
%the concept of focusing sequences, developed
%in \cite{KLS3} for the Maxwell's equations for the Dirac equation.
\smallskip

\noindent {\bf  7.1\ Controlling supports of  generalized waves.}
% We recall that by Lemma \ref{lem: equivalence of data},
%$C_0(D)$ and the induced bundle structure on $\p M$
%determine $\Lambda_\Gamma$.
%Because of this we analyze the local boundary condition
%given by $P_\Gamma$.
We start by fixing certain notations.  
Let  
%the times $T_0$ and $T_1$ 
%satisfy
\beq
\label{T} 
 T_1 =T+4\sup\,\{ {\rm rad}\,(M,\Sigma)\ :\ \Sigma\not=\emptyset\}.
\eeq 
Then $T_1$ is determined by the boundary data
and  $T_1-T> \hbox{diam}\,(M)$.

Let $\Sigma_j\subset\partial M$ be non-empty open 
disjoint sets, $1\leq j\leq J$ and $\tau_j^-$ and 
$\tau_j^+$ be positive times with 
\[ 
 0<\tau_j^-<\tau_j^+\leq T_1-T, 
 \quad 1\leq j\leq J. 
\] 
Let  $S = S(\{\Sigma_j,\tau_j^-, 
\tau_j^+\}_{j=1}^J) \subset M$ be the intersection 
of slices, 
\vspace{-5mm}

\begin{equation}\label{def of S} 
  S =\bigcap_{j=1}^J  \left( M(\Sigma_j,\tau_j^+)\setminus 
M(\Sigma_j,\tau_j^-) \right). 
\end{equation} 
Our first goal is to find, using 
the boundary data, 
whether the set $S$ contains an open ball or not. 
To this end, we give the following definition.

\begin{definition}\label{support sources} 
{\rm The set $Z = Z(\{\Sigma_j,\tau_j^-, 
\tau_j^+\}_{j=1}^J)$ is the set of the  generalized sources 
$\hat f\in \F^\infty$ that generate  the
waves ${\mathcal W}{\hf}(t)$ 
 with
\begin{enumerate} 
\item $\supp({\mathcal W}{\hf}(T_1))\subset M(\Sigma_j,\tau_j^+)$ for 
all $j$, $1\leq j\leq J$, 
\item ${\mathcal W}{\hf}(T_1)=0$ in $M(\Sigma_j,\tau_j^-)$ 
for 
all $j$, $1\leq j\leq J$. 
\end{enumerate} 
}
\end{definition}

\begin{theorem}\label{is in Z th} 
Given the
induced bundle structure
on $\p M$
 and  the  response operator ${\Lambda}_{\Gamma}$, 
we can 
determine whether a given boundary source $\hat f 
\in \F^\infty$ is in 
 $Z$ or not. 
\end{theorem}

\noindent {\em Proof:} Let 
$\hat f=(f_k)_{k=0}^\infty \in \F^\infty$
 be a generalized source. 
Consider first the question whether $\supp({\mathcal W}{\hat f}(T_1))\subset
M(\Sigma_j,\tau_j^+)$ or equivalently,
$\supp(\Pi_+ {\mathcal W}{\hat f}(T_1))\subset
M(\Sigma_j,\tau_j^+)$
and $\supp(\Pi_- {\mathcal W}{\hat f}(T_1))\subset
M(\Sigma_j,\tau_j^+)$.
 By Theorem \ref{local control th}
this is equivalent to the existence of generalized sources,
\[
 \hat h^+ =(h_\ell^+)_{\ell=0}^\infty\hbox{ and }
 \hat h^- =(h_\ell^-)_{\ell=0}^\infty
,\quad 
 h_\ell^\pm \in C^\infty_0(\Sigma_j \times (0,\tau_j^+),\, P_{\Gamma} W), 
\]
such that 
\beq 
\label{20.11.1}
 \lim_{k,\ell \to \infty} \| \Pi_+ (u^{f_k}(T_1) 
-u_t^{h_\ell^+}(\tau_j^+))\|_{L^2(M,V)} =0,\\ \label{20.11.1b} 
 \lim_{k,\ell \to \infty} \| \Pi_- (u^{f_k}(T_1) 
-u_t^{h_\ell^-}(\tau_j^+))\|_{L^2(M,V)} =0. 
\eeq 
By Theorem \ref{blacho}, we can evaluate the left-hand sides of
(\ref{20.11.1}) and (\ref{20.11.1b})  using 
the induced bundle structure
(\ref{b-structure}) on $\p M$ and 
the operator $\Lambda_{\Gamma}$.
Thus, for a given $\hat f$, we can
 verify if  condition {\it 1.}
is satisfied.

Next, consider condition {\it 2}. As
\[ 
 (i\partial_t+D){\mathcal W}{\hf}=0 \hbox{ in} \,\, M\times\R,  \quad
 P_{\Gamma} \left({\mathcal W}{\hf}|_{\partial M\times \R}\right)=0,  
\] 
condition {\it 2.} implies, due to the finite propagation speed, that
\[ 
{\mathcal W}{\hf}=0\mbox{ in $K_j=\{(x,t)\in M\times\R:\quad \dist(x,\Sigma_j) 
+|t-T_1|<\tau_j^-\}$}, 
\] 
for all $j=1,\ldots,J$. 
{\newwtextt In particular,  
\beq
\label{66b}
 {\mathcal W}{\hat f}|_{\Sigma_j\times (T_1- \tau_j^-, \, T_1+ \tau_j^-)}=
 \lim_{j\to \infty}\Lambda_{\Gamma} f_j|_{\Sigma_j\times (T_1- \tau_j^-, \, T_1+ \tau_j^-)} =0.
\eeq
 On the other hand,
if   (\ref{66b})  takes place,
then Theorem \ref{UCP} implies
that ${\mathcal W}{\hf}=0$  in $K_j$. As
 condition (\ref{66b}) 
can be verified given the induced bundle structure
 on $\p M$ and  the response operator $\Lambda_{\Gamma}$, 
 we can verify condition {\it 2.} }
%. Thus $\hat f=(f_j)_{j=1}^\infty$
% satisfies condition {\it 2}
%if $\lim_{j\to \infty}(I- P_{\Gamma})f_j=0$ 
%in   
%$\Sigma_j\times (T_1- \tau_j^-, \, T_1+ \tau_j^-)$
%for all $j=1,\dots,J$.
\hfill$\Box$

\begin{theorem}\label{alternative} 
Let $S$ and $Z$ be as above. The following 
alternative holds: 
\begin{itemize} 
\item [1.] If $S$ contains an open ball, then 
${\rm dim}(Z)=\infty$, 
\item [2.] If $S$ does not contain an open ball, then 
$Z = \{0\}$. 
\end{itemize} 
\end{theorem}

\noindent {\em Proof:} {\it 1.}
Let  $B\subset S$ be an open ball and  
 $0\neq\varphi\in C^\infty_0( B,V)$. By  the
global controllability, Theorem \ref{global 
control th}, there is
 $\hat f \in {\F}^{\infty}$ 
such that 
\ba
%\label{20.11.2} 
 {\mathcal W}{\hat f}(T_1)=\varphi.
\ea
As  $C^\infty_0( B,V) \subset {\cal D}(D_{\Gamma}^{\infty})$, then 
$\hat f\in Z$. As $\phi\in C_0^\infty(B,V)\setminus \{0\}$ is arbitrary, this implies  that
  $\hbox{dim}(Z) = \infty$.

{\it 2.} Assume $S$ does 
not contain an open ball and  suppose that there is  
 $0 \neq \hat f\in Z$. 
Then, by conditions {\it 1., 2.} in  Definition 
\ref{support sources}, 
\ba
& & {\rm supp}({\mathcal W}{\hat f}(T_1))\subset 
\bigcap_{j=1}^J M(\Sigma_j,\tau_j^+) = S^+, \\
& &{\mathcal W}{\hat f}(x,T_1)=0,\quad \mbox{for $x\in\bigcup_{j=1}^J 
M(\Sigma_j,\tau_j^-)=S^-$.} 
\ea
Thus ${\rm supp} ({\mathcal W}{\hat f}(T_1)) \subset S^+\setminus S^-$.
As $S$ does not contain a ball, then
$S^+\setminus S^-$ is nowhere dense. 
Since ${\mathcal W}{\hat f}(T_1)$ is smooth,  it 
 vanishes  everywhere in $M$. This is a contradiction with
the assumption that $\hat f\not=0$ in $\F$.
\hfill$\Box$  
\smallskip

We are now ready to construct the set of the boundary 
distance functions. For each $x\in M$, 
the corresponding boundary 
distance function, $r_x \in C(\p M)$ is  given by 
\[ 
 r_x: \p M\to\R_+,\quad r_x(z)=\dist(x,z), \quad 
z \in \partial M.
\] 
%In fact, $r_x\in \hbox{Lip}\,(\p M)$ with the Lipschitz constant equal to one.
The boundary 
distance functions define {\it the 
 boundary distance map} ${\mathcal R}:M\to C(\p M)$,
${\mathcal R}(x)=r_x$, which is continuous and injective
(see \cite {Ku5,KKL}).  Denote by
\[ 
 {\mathcal R}(M)=\{r_x\in C(\partial M): \, x\in M\} 
\] 
the image
of ${\mathcal R}$. It is known (see \cite{Ku5,KKL})
  that, given the set
 ${\mathcal R}(M) 
\subset C(\partial M)$, 
we can endow it, in a constructive way, 
 with  
a  differentiable structure and a metric tensor $\tilde g$, 
so that $({\mathcal R}(M),\tilde g)$ becomes a manifold that is
 isometric
to $(M,g)$, 
\[ 
 ({\mathcal R}(M),\tilde g)\cong (M,g). 
\] 
We complete the reconstruction of $(M,g)$ with the  following result:

\begin{theorem} 
Let  the
induced bundle structure
(\ref{b-structure}) on $\p M$
and  the response operator ${\Lambda}_{\Gamma}$
be given. Then, for 
any $h\in C(\partial M)$, we can find whether 
$h\in{\mathcal R}(M)$ or not. 
\end{theorem}

\noindent {\em Proof:} % The proof is identical to that of \cite{KLS3}:
Consider a function $h\in C(\p M)$.
We have that $h\in {\mathcal R}(M)$ if and only if 
there is $x\in M$ such that $h=r_x$. 
Observe that if such $x$ exists then, for any $\e>0$, we have
$h(x)-\e<d(y,z)<h(x)+\e$ for any $z\in \p M$ and
$y\in B(x,\e)$, where $B(x,\e)\subset M$ is $x$-centered 
ball of radius $\e$. Therefore,
for any $J \in \Z_+$,
any set of points $z_j\in \p M$, and non-negative numbers $\tau_j^-,\tau_j^+$ satisfying
$\tau_j^-<h(z_j)<\tau_j^+$, $j=1,\dots,J,$ the set 
\beq\label{S-ehto}
S=S(\{\Sigma_{j},\tau_{j}^-, 
\tau_{j}^+\}_{j=1}^{J})\subset M \quad \hbox{contains an open ball,}
\eeq
if  $\Sigma_j\subset \p M$ 
are sufficiently small neighborhoods of
of $z_j$.

%  Indeed, let
%\ba
%\e = \min_j(|h(z_j)-\tau_j^{\pm}|)
%\ea
%Then, for any $\Sigma_j$ with $\diam_{\p M}(\Sigma_j) < \e/2,$
%\ba
%B_{\e/2}(x) \subset S.
%\ea
On the other hand, assume that  for all  $J$,
$z_j\in \p M$, $\tau_j^-,\tau_j^+$ satisfying
$\tau_j^-<h(z_j)<\tau_j^+$, and sufficiently small $\Sigma_j$, 
 condition (\ref{S-ehto}) is satisfied.
%for all sufficiently small $\Sigma_j$
Then there is a point $x\in M$ such that
\ba
|d(x,z_j)-h(z_j)|\leq \tau_j^+-\tau_j^-+\diam(\Sigma_j),\quad j=1,2,\dots, J.
\ea
Thus, letting $\tau_j^-\to h(z_j)$, $\tau_j^+\to h(z_j)$ and,
$\diam(\Sigma_j)\to 0$, we see that
there is $x_J\in M$ such that $d(x_J,z_j)=h(z_j)$ for all
$j=1,2,\dots, J$.

Denote by $\{z_j\}_{j=1}^\infty$ a countable 
dense set on $\p M$, and let $x_J$, $J\in \Z_+$ be the above defined points
corresponding to the first $J$ 
points  $z_j, \, j=1,\dots, J,$ of this set.  This gives us a 
sequence $(x_J)_{J=1}^\infty$. 
As $M$ is compact, this sequence has a limit point
 $x \in M$. By density of $\{z_j\}_{j=1}^\infty$, $h=r_x$.
 % Indeed, if a point $x$ exists, the
%set $S$ contains an open ball centered at $x$. If no such a point
%$x$ exists, we see using closedness of
%$ {\mathcal R}(M)$ in the topology of pointwise convergence
%of $C(\p M)$ we see that
%there are $z_j$ and $\tau_j^\pm$ such that $S$ is empty.

By Theorem \ref{alternative},
condition  (\ref{S-ehto}) is equivalent to  the condition
%\beq\label{ehto 1}
\ba
\dim(Z(\{\Sigma_{j},\tau_{j}^-, 
\tau_{j}^+\}_{j=1}^{J}))\not =0.
\ea
%\eeq
By Theorem \ref{is in Z th}, this  can be verified using the boundary data.
  \hfill$\Box$ 

 \smallskip

Combining the above theorem  with Lemma \ref{lem: equivalence of data},
 we obtain the following result:

\begin{corollary} 
The 
induced bundle structure
(\ref{b-structure}) on $\p M$ and
 the Cauchy data set $C_0(D)$
determine the Riemannian manifold $(M,g)$ 
 up to an isometry.
\end{corollary}

\section{Focusing sequences. 
Proof of Theorem \ref{main1}.}

\noindent {\bf  8.1\ Waves and delta-distributions.}
In the previous section it was shown 
 that, using 
the induced bundle structure
 on $\p M$ and
the response operator $\Lambda_{\Gamma}$,
 we can control the supports
 of the  
waves ${\mathcal W}{\hat{f}}(t)$. In this section the goal is 
to construct a sequence 
of  sources, 
$(\hat{f}_l)_{l=1}^\infty \subset \F^{\infty}$ 
such that,
when $l \to \infty$, the corresponding waves ${\mathcal W}{\hat{f}_l}(t)$
concentrate, at $t=T_1$,
where $T_1$ is defined by (\ref{T}), at a single 
point $y \in M^{\rm int}$. 
These waves
% For $t \geq T_1$,  
%these waves behave like 
%the waves generated by point sources at $y$, a fact 
%that turns out
turn to be useful for reconstructing 
the bundle $V$ and the Dirac-type operator $D$.

In the following, let $\ud_y$ denote the Dirac delta-distribution
at $y\in M^{{\rm int}}$, i.e., 
\[ 
 \int_M\ud_y(x)\phi(x)dV_g(x)=\phi(y),\quad \hbox{for }\phi\in C^\infty_0(M). 
 \]

Since the Riemannian manifold $(M,g)$ is already found, 
for any $y \in M^{{\rm int}}$ we can choose a sequence
 $\Sigma_{jl}\subset\partial M$, 
$\,0<\tau_{jl}^-<\tau_{jl}^+\leq T_1-T$, so that,
for $S_l=S(\{\Sigma_{jl},\tau_{jl}^+,\tau_{jl}^-\}_{j=1}^{J(l)})$,
\beq
\label{20.11.3} 
 S_{l+1}\subset S_l,\quad \bigcap_{l=1}^\infty S_l 
 = \{y\}. 
\eeq 
Let
 $Z_l=Z(\{\Sigma_{jl},\tau_{jl}^-, 
\tau_{jl}^+\}_{j=1}^{J(l)})$ be  
the corresponding set 
of the generalized sources defined in Definition 
\ref{support sources}.

\begin{definition} {\rm 
Let $S_l, \, l=1,2,\dots,$ satisfy (\ref{20.11.3}).   
 We call the sequence $ (\hat f_l)_{l=1}^\infty,$
 with 
$\hat f_l\in Z_l$, a {\it focusing sequence} 
 of generalized sources of the order  $s \in \R_+$  
(for brevity, a focusing sequence),  if 
there is a distribution-valued section $0 \not=A_y \in {\mathcal D}'(M,\,V)$ 
supported at the  point $y\in M^{int}$ such that 
\[ 
 \lim_{l\to \infty} \bbra {\mathcal W}{\hat f_l}(T_1) 
,\eta\ccet  = \bbra A_y,\eta\ccet,\quad \hbox{for all $\eta\in {\mathcal D}(|D_{\Gamma}|^s)$.} 
\] 
}
\end{definition}

\begin{lemma} 
Assume that the induced bundle structure
(\ref{b-structure}) on $\p M$ and  the response operator 
$\Lambda_{\Gamma}$ be given.
Then, for any $s \in \Bbb{R}_+$ and any sequence  
of generalized sources,
$(\hat f_l)_{l=1}^\infty\subset \F^\infty$, we can determine if
$(\hat f_l)_{l=1}^\infty$ is a focusing sequence of the order $s$ or not.  
\end{lemma}

\noindent {\em Proof:} By Theorem \ref{global control th},  any  
$\eta\in {\mathcal D}(|D_{\Gamma}|^s)$
can be represented as $\eta={\mathcal W}{\hat h}(T_1),$
 $\hat h\in\F^s$. Thus, 
if $(\hat f_l)_{l=1}^\infty$ is a focusing sequence 
of the order $s$, then
 the limit, 
\begin{equation}\label{limit1} 
 \bbra A_y,\, {\mathcal W}{\hat h}(T_1)\ccet=\lim_{l\to \infty} \bbra {\mathcal W}{  
\hat f_l}(T_1),\, {\mathcal W}{\hat h}(T_1)\ccet, 
\end{equation} 
exists for all $\hat h\in \F^s$. By
Theorem \ref{blacho},  the existence of this
limit can be 
verified if we are given the induced bundle structure
on $\p M$ and the operator
 $\Lambda_{\Gamma}$.

Conversely, assume that {\newwtext 
we have a sequence 
$ (\hat f_l)_{l=1}^\infty,$
 with 
$\hat f_l\in Z_l$ such that
the limit (\ref{limit1})  
exists for all $\hat h\in\F^s$ 
and that this limit is non-zero for some  $\hat h\in\F^s$.} 
Then, by the principle of the uniform boundedness,  
the mappings 
\[ 
 \eta \mapsto \bbra {\mathcal W}{ \hat f_l}(T_1) 
,\, \eta\ccet, \quad l \in \Bbb{Z}_+, 
\] 
form a uniformly bounded family in the dual of 
${\mathcal D}(|D_{\Gamma}|^s)$. By the
Banach-Alaoglu theorem, we find   a weak$^*$-convergent 
subsequence in this space with the limit
\[ 
A_y:=\lim_{l\to \infty} {\mathcal W}{ \hat f_l}(T_1)
\in \bigg({\mathcal D}(|D_{\Gamma}|^s) 
\bigg)'. 
\] 
This  gives us  the  required non-zero distribution-valued
section in ${\mathcal D}'(M, V)$.
As $\hat f_l\in Z_l$, then $\hbox{supp}(A_y) = \{y\}$.
\hfill$\Box$ \smallskip

Next we consider the  representation of 
$A_y$ in a local trivialisation near $y$.
Since ${\rm supp}(A_y)=\{y\}$, 
$A_y$ is a finite linear combination,
with coefficients from $\pi^{-1}(y)$,
 of the Dirac delta-distribution at $y$ and its derivatives. The 
role of the smoothness index $s$ is just to select the 
order of this distribution, as is seen in the 
following result.

\begin{lemma} 
\label{Lm2.15}
Let $ A_y$, $y\in M^{\rm int}$ be a limit of a focusing sequence of the order $s,$ 
$\frac n2<s<\frac {n+1}2$.
Then $A_y$ is of the form 
\ba%\label{Ay} 
 A_y(x)= \lambda\ud_y, 
\ea
where $\lambda\in \pi^{-1}(y)$. 
Furthermore, for any $ \lambda \in \pi^{-1}(y)$ there is a 
focusing sequence
$(\hat f_l)_{l=1}^\infty$ {\newwtext such that the
waves ${\mathcal W} {\hat f}_l(T_1) \rightarrow \lambda\ud_y$ 
 in ${\mathcal D}'(M,\,V)$} when $l \to \infty$. 
\end{lemma} 

 \noindent {\em Proof:}  
 As for $s \geq 0$
\ba
H^s_{\hbox{\tiny comp}}(M,\,V) \subset {\mathcal D}(|D_{\Gamma}|^s) 
\subset H^s(M,\,V),
\ea
the embedding theorems imply that, for $\frac{n}{2} <s< \frac{n+1}{2}$,
the delta-distributions $\lambda\ud_y, \, \lambda \in \pi^{-1}(y)$, 
are in $\left({\mathcal D}(|D_{\Gamma}|^s)\right)'$ while the
derivatives of the delta-distribution are
not in $\left({\mathcal D}(|D_{\Gamma}|^s)\right)'$.

To prove  
 the assertion we have only to verify that a  desired 
focusing sequence 
exists
for  any $\lambda \in \pi^{-1}(y)$. This can be done by an approximation of $\lambda \ud_y$
by smooth sections with supports converging to $y$ 
and an application of  the global controllability Theorem 
\ref{global control th}.
\hfill$\Box$ 

\smallskip

\noindent {\bf  8.2. Construction of the bundle structure.}

\begin{lemma} \label{gauge2}
Assume that 
we are given the induced bundle structure
(\ref{b-structure}) on $\p M$ and the response operator $\Lambda_{\Gamma}$.
Then, for any 
$y_0 \in M^{{\rm int}}$, we can find integers 
$k(1),\dots,k(d)$ and a neighborhood $U\subset M^{\rm int}$ of $y_0$
such that the eigenfunctions
$\phi_{k(1)}(y),\dots,\phi_{k(d)}(y)$ of $D_{\Gamma}$ 
form a basis in $\pi^{-1}(y)$, $y\in U$.
Moreover, for any ${\hat h} \in \F^{\infty}$, 
we can find the coefficients $\alpha_l(y,t)$ 
of the decomposition
\beq
\label{25.11.6}
{\mathcal W}\hh(y,t) = \sum_{l=1}^d \alpha_l(y,t) \phi_{k(l)}(y), \quad y \in U,
\ t\in \R.
\eeq
\end{lemma}

\noindent {\em Proof:} Let $U_1 \subset M^{{\rm int}}$ be a sufficiently small neighborhood 
of a point $y_0\in M^{\rm int}$ so that there is a trivialisation 
 $\Phi:
\pi^{-1}(U_1)\to U_1\times \C^d$.
For all $y\in U$ we choose  $d$
families of focusing sequences,
$(\hat f_k^{j,y})_{k=1}^\infty$,  $ j=1,\dots,d$, 
such that ${\mathcal W}\hat f_k^{j,y}(T)$ converge to
some distributions $\lambda_y^j \ud_y,\,$
 $\lambda_y^j\in \pi^{-1}(y)$, $y\in U_1$.
 
 {\newwtextt Using the operator $\Lambda_{\Gamma}$
and the induced bundle structure
 on $\p M$ we can evaluate, by Theorem \ref{blacho},
the inner products
\beq\label{eq: Matti A}\quad\quad
\lim_{k \to \infty}\bra \bra  {\mathcal W}\hh(T),  {\mathcal W}\hf_k^{j,y}(T) \cet \cet 
=
\bra\bra {\mathcal W}\hh(T), \lambda_y^j\ud_y\cet\cet
=
\bra {\mathcal W}\hh(y,T),\lambda_y^j\cet_y,\hspace{-1cm}
\eeq
where $j=1,\dots,d.$
The  set $\{ {\mathcal W}\hh(T):\ \hh\in \F^\infty\}$
contains $C^{\infty}_0(M, V)$. Therefore,
$\{ {\mathcal W}\hh(y, T):\ \hh\in \F^\infty\}= \pi^{-1}(y)$ and 
it is possible to verify, using (\ref{eq: Matti A}), whether the set
$\{\la_y^j, \, j=1, \dots, d\}$ is a basis in $\pi^{-1}(y)$ for
$y$ lying in some neighborhood $U \subset U_1$ of $y_0$. }
%Thus, in
%particular,  $ \{{\mathcal W}\hh(y, T): \,\hh\in \F^\infty\} = \pi^{-1}(y)$.
%Thus we can verify, 
%whether
%for given families of the focusing
%sequences
% all inner products (\ref{eq: Matti A}), $\hat h\in \F^\infty$
%are $C^\infty$-smooth functions on $M^{\rm int}$, or equivalently,
%  the set
%$\{\la_y^j, \, j=1, \dots, d\}$ is linearly independent
%for all $y$ in some sufficiently small
% neighborhood $U_1$ of $y_0$.
Moreover, we can verify whether $\la_y^j$, considered as 
functions of $y$, define $C^\infty$-smooth sections in
$U$.

On the other hand, Lemma \ref{Lm2.15}
guarantees the existence of  focusing sequences 
 $(\hat f_l^{j,y})_{l=1}^\infty$ such that the corresponding
$\la_y^j$ are linearly independent and smooth in $U$.
%In the sequel, let $\lambda_y^j$, $j=1,\dots,d$ be some 
%smooth sections in $C^\infty(U_1,V)$ that  
%form a basis in $\pi^{-1}(y)$ for any $ y \in U_1$.

For any  
$\hh\in \F^\infty$, $t\in \R$, consider the map $K_t:\F^\infty
\to C^\infty(U,\C^d)$,
\ba
%\label{map}
K_{t}\hh(y)= (\bra {\mathcal W}\hh(y,t),\lambda_y^1\cet_y), \dots,
\bra {\mathcal W}\hh(y,t),\lambda_y^d\cet_y),
\quad y\in U.\hspace{-1cm}
\ea
%Note that $K_t\hh(y)$, $y\in U_1$
%can be evaluated using the induced bundle structure
%on $\p M$ and the response operator $\Lambda_{\Gamma}$.
By Proposition
\ref{5.1}, we can find the generalized sources $\hh_k$, $k=1,2,\dots$
such that 
${\mathcal W}{\hh_k}(T)= \phi_k$.
Thus, for any  $d$ numbers $k(1),\dots,k(d)\in \Z_+$,
we can evaluate the matrix
\ba
%\label{25.11.4}\quad\quad
{\mathcal E}(y):=[e_{j,l}(y)]_{j,l=1}^d, \quad 
e_{j,l}(y) = \bra \phi_{k(l)}(y),\, \lambda^j_y \cet_y,\quad y\in U_1.
\ea
As for any $y_0 \in M^{{\rm int}}$ the vectors 
$\{ \phi_k(y_0)\}_{k=1}^{\infty}$ span 
$\pi^{-1}(y_0)$, it is possible to select $k(1),\dots,k(d)$
so that the system $\Phi$ of $d$ 
eigenfunctions,
$\Phi=\left(\phi_{k(1)}, \dots,\phi_{k(d)}\right)$ 
is linearly independent at $y_0$. As ${\mathcal E}(y)$
is smooth, we can choose  $U\subset U_1$ to be a neighborhood 
of $y_0$ where ${\mathcal E}(y)$ is invertible.

Let $\hh\in \F^\infty$. Then the coefficients $\alpha_i(y,t)$ 
in (\ref{25.11.6}) can be found by computing
the vector
$(\alpha_1(y,t),\dots, \alpha_d(y,t))={\mathcal E}(y)^{-1}K_t\hh(y)$.
\proofbox\smallskip

Next we consider an eigenfunction basis near the boundary.
Due to the boundary condition
$P_{\Gamma} ( \phi_k|_{\p M})=0$, the eigenfunctions $\phi_k$ 
do not 
span $\pi^{-1}(y_0),$ $y_0 \in \p M$.
Because of this, we will use near $\p M$  a combination of the 
eigenfunctions $\phi_k^+$ of $D_{\Gamma,+}$ 
and $\phi_j^-$ of $D_{\Gamma,-}$, where 
 $D_{\Gamma,+}=D_{\Gamma}$ and
$D_{\Gamma,-}=D_P$ with $P= (I-P_{\Gamma})$ 
(compare with \cite{AK2LT} where 
a coordinate system near $\p M$ 
was constructed using a combination of the Dirichlet 
and Neumann eigenfunctions). 
We denote the spaces of the generalized sources
corresponding to the operators $D_{\Gamma,\pm}$
by $\F_\pm$ and $\F^s_\pm$ and by ${\mathcal W}_\pm$
be propagators (\ref{gener}) for $D_{P,\pm}$.
We note that our previous results concerning the unique 
continuation, controllability, and the construction of
focusing sequences can be proven using the boundary 
condition $I-P_{\Gamma}=\frac 12(I-\Gamma)$ instead
of $P_{\Gamma}=\frac 12(I+\Gamma)$.

To simultaneously consider the
 operators $D_{\Gamma,+}$ and  $D_{\Gamma,-}$, 
we need the following modification of Theorem \ref{blacho}:

\begin{lemma}
\label{blacho2} Assume that we are given 
the induced bundle structure on $\p M$ and
the response operator $\Lambda_{\Gamma}$. Then for
$
f_+ \in C^{\infty}_0(\p M \times \R_+, \, P_{\Gamma} W)$ and
$f_- \in C^{\infty}_0(\p M \times \R_+, \, (1-P_{\Gamma}) W)$
we can evaluate the inner products 
\ba
\bbra u^{f_+}(t),\, u^{f_-}(s) \ccet,\quad
\hbox{for }t,s \geq 0. 
\ea
Here $u^{f_+}$ is the solution of (\ref{IBVP0})-(\ref{IBVP0line 2}) with
$P=P_{\Gamma}$ and $u^{f_-}$  is the solution of (\ref{IBVP0})-(\ref{IBVP0line 2}) with
$P=I-P_{\Gamma}$. 
\end{lemma}

\noindent {\em Proof:}
Let
$
I(s,t) = \bbra u^{f_+}(t),\, u^{f_-}(s) \ccet.
$
Similarly to the proof of Theorem \ref{blacho}, 
we use integration by parts to see that 
\beq\label{ODY 1st}\quad\quad
(\p_s+\p_t)I(s,t)=i\bbra  Du^{f_+}(t),u^{f_-}(s)\ccet
-i\bbra  u^{f_+}(t),Du^{f_-}(s)\ccet  
\eeq
\ba
& & = i \int_{\p M} 
\bra \gamma(N)\Lambda_{\Gamma}^+  f_+(t),\, \Lambda_{\Gamma}^-  f_-(s) \cet\, dA_g.
\ea
Here $\Lambda_{\Gamma}^{\pm}$ are the response operators for 
$D_{\Gamma,\pm}$. By Lemma \ref{lem: equivalence of data}
the operator $ \Lambda_{\Gamma}= \Lambda_{\Gamma}^+$ and the
induced bundle structure determine $ \Lambda_{\Gamma}^-$.
As $I(0,t)=I(s,0)=0$, we can determine $I(s,t)$ by solving
the differential equation (\ref{ODY 1st}) along the characteristics
$t-s=c$.
\hfill$\Box$ \smallskip

Let us denote by $\psi_k$ of the eigenfunctions
of $D_{\Gamma,+}$ and $D_{\Gamma,-}$, renumbered 
so that $\psi_{2k-1}=\phi^+_{k}$ and  $\psi_{2k}=\phi^-_{k}$, $k\geq 1$. 
To generalize
Lemma \ref{gauge2}, we need the following result:

\begin{lemma}
\label{C infinity} ${\mathcal D}(D_{\Gamma,+}^\infty)+ 
{\mathcal D}(D_{\Gamma,-}^\infty)=
C^\infty(M,V)$.
 \end{lemma}

\noindent {\em Proof:} 
For  $z\in \p M$ let be a sufficiently small
 ball $B=B(z,r)\subset M$
%Let $z\in \p M$ and $r>0$ be such
%that in the ball $B=B(z,r)\subset M$
such that the boundary normal coordinates
 and
 a local trivialization
$\Phi_B:
\pi^{-1}(B)\to B\times \C^d$ are well defined in $B$.
It suffices to show that  any $u\in C^\infty(M,V)$ with $\hbox{supp}(u) \subset
B(z,r/2)$ can be written as
\beq\label{eq: u plusminus}
u=u_++u_-,\quad \hbox{where }
u_+\in {\mathcal D}(D_{\Gamma,+}^\infty), \ \ 
u_-\in {\mathcal D}(D_{\Gamma,-}^\infty).
 \eeq

To prove (\ref{eq: u plusminus}), 
%we first consider
%a consequence of Seeley's extension theorem.
%To this end, 
let $u^j\in C^\infty(\p M,W)$, $j=0,1,2,\dots$ be supported
on $B(z,r/2)\cap \p M$.
By using Seeley's extension theorem \cite{seeley} in $B$,
%the boundary normal coordinates, 
there is $u\in C^\infty(M,V)$ supported on $B(z,r/2)$
such that
\ba
\p_N^ j u|_{\p M}=u^j.
\ea
Moreover, using representations (\ref{dirac 1}), (\ref{dirac-type, maximal})
of the Dirac operator in the local trivialization $\Phi_B$,
we see that
\beq\label{recurrent dirac}
D^ju |_{\p M}=\gamma(N) u^j +A_j(u^0,u^1,\dots,u^{j-1}),
\eeq
where $A_j$ is a differential operator.
Let $w^j\in C^\infty_0(\p M,W)$, $j=0,1,2,\dots,$ be supported
in $B(z,r/2)\cap \p M$.
Solving recurrently equations (\ref{recurrent dirac})
we see that there is $u\in C^\infty(M,V)$ such that
$D^ju |_{\p M}=w_j$.

Let $v\in C^\infty(M,V)$. Then
 $v\in  {\mathcal D}(D_{\Gamma,\pm}^\infty)$
if and only if 
\beq\label{pm domain}
P_{\Gamma,\pm} (D^jv)|_{\p M}=0,\quad j=0,1,2,\dots.
\eeq 
Thus, if
$u\in C^\infty(M,V)$ supported in
$B(z,r/2)$
is of the form (\ref{eq: u plusminus}),
we have by (\ref{pm domain}) that
\ba 
& &D^ju_+|_{\p M}=(1-P_{\Gamma})( D^ju|_{\p M})=:w_j^+,\\
& &D^ju_-|_{\p M}=P_{\Gamma} (D^ju|_{\p M})=:w_j^-.
\ea
On other hand, let $u\in C^\infty(M,V)$ and
 $u_+,u_-\in C^\infty(M,V)$ be functions
such that $ D^ju_\pm |_{\p M}=w_j^\pm$.
Then $ u_\pm\in {\mathcal D}(D_{\Gamma,\pm}^\infty) $
and (\ref{eq: u plusminus}) is satisfied.
\hfill$\Box$

\begin{lemma}
\label{gauge3}
Assume that 
we are given the induced bundle structure
(\ref{b-structure}) on $\p M$ and the response operator 
 $\Lambda_{\Gamma}$.
Then, for any 
$y_0 \in M$, we can find integers 
$k(1),\dots,k(d)$ and a neighborhood $U\subset M$ of $y_0$
such that the eigenfunctions $\psi_{k(1)}(y),\dots,\psi_{k(d)}(y)$ 
of 
  $D_{\Gamma,+}$ or $D_{\Gamma,-}$ form a basis of $\pi^{-1}(y)$, $y\in U$.
Moreover, for any ${\hat h} \in \F^{\infty}_\pm$, 
we can find coefficients $\alpha_l^{\pm} (y,t)$ 
such that 
\beq
\label{27.11.2}
{\mathcal W}_\pm {\hh_\pm}(y,t) = \sum_{l=1}^d \alpha_l^{\pm}(y,t) \psi_{k(l)}(y), \quad y \in U,\ t\in \R.
\eeq
\end{lemma}

\noindent {\em Proof:}
First we show that,
for any $y_0 \in \p M$, the space $\pi^{-1}(y_0)$ is 
spanned by $\{\phi_k^+ (y_0)\}_{k=1}^{\infty}\cup 
\{\phi_k^- (y_0)\}_{k=1}^{\infty}.$
To  this end, let $u\in C^\infty(M,V)$. By Lemma \ref{C infinity} there are
$u_+\in {\mathcal D}(D_{\Gamma,+}^\infty)$ and 
$u_-\in {\mathcal D}(D_{\Gamma,-}^\infty)$ such that
$u=u_++u_-$. Consider the eigenfunction expansions
\beq\label{eq: eigen expansion}
u_\pm=\sum_{k=1}^\infty a^\pm_j \phi_k^\pm
\eeq
that converge in  $ {\mathcal D}(D_{\Gamma,\pm}^\infty)$.
By elliptic regularity in Theorem \ref{self-adjointenss},
 the series (\ref{eq: eigen expansion})
converges in  $C^\infty(M,V)$. In particular, any
$\mu\in \pi^{-1}(y_0)$ can be written in the form
$\mu=\sum_{k=1}^\infty a^+_j \phi_k^+(y_0)+ 
\sum_{k=1}^\infty a^-_j \phi_k^-(y_0).$

Now, let $y_0 \in \p M$. 
{\newwtextt Arguing as in the proof of Lemma \ref{gauge2} but using Lemma \ref{blacho2}
instead of Theorem \ref{blacho}, we find $k(1),\dots,k(d)\in \Z_+$
and neighborhoods $U$ of $y_0$ such that
$\psi_{k(1)}(y),\dots,\psi_{k(d)}(y)$ form a smooth basis in $\pi^{-1}(U\cap M^{\rm int})$.
Moreover, we can find, using the boundary data, the coefficients $\alpha_l^{\pm} (y,t)$
in (\ref{27.11.2}) for any $y \in U\cap M^{\rm int}$.  Then, by Lemma \ref{C infinity}, 
$\psi_{k(1)}(y),\dots,\psi_{k(d)}(y)$ form a smooth basis in $\pi^{-1}(U)$
if and only if the coefficients $\alpha_l^{\pm} (y,t)$ may be continued smoothly
on the whole $U$ for any ${\hat h} \in \F^{\infty}_\pm$.  This can be verified using the boundary data.
}
%Using Lemma \ref{blacho2}, we can find exactly
%as in the proof of Lemma \ref{gauge2}  all those $k(1),\dots,k(d)\in \Z_+$
%and a neighborhood $U$ of $y_0$ such that
%$\psi_{k(1)}(y),\dots,\psi_{k(d)}(y)$ are linearly independent
%for $y\in U\cap M^{\rm int}$. Now
%$\psi_{k(1)}(y),\dots,\psi_{k(d)}(y)$ are linearly
%independent for  all $y  \in U\cap \p M$ if and only
%if for any $w={\mathcal W}_+\hat h_+(t)+{\mathcal W}_-\hat h_-(t)$,
%$\hat h_\pm\in \F_\pm^\infty,$ $t\in \R$ 
%the coefficients in representation 
%(\ref{27.11.2}) of $w$ can be extended to smooth functions up to the boundary
%$U\cap \p M$.
%This can be verified using the boundary data, and thus we can
%find all those $k(1),\dots,k(d)\in \Z_+$ and $U$ for which the assertion is valid.
\proofbox\smallskip

The system  $\Phi=(\psi_{k(1)}(y),\dots,\psi_{k(d)}(y))$ of eigenfunctions
and the neighborhood $U$ in Lemma \ref{gauge3} define
a local trivialization of the vector bundle $V$, namely
the map $A_{\Phi,U}:\pi^{-1}(U)\to U\times \C^d$,
\ba
A_{\Phi,U}(v)=(\alpha_j(y))_{j=1}^d,\ \ v\in \pi^{-1}(y),\hbox{ when }
v= \sum_{l=1}^d \alpha_l(y) \psi_{k(l)}(y).\hspace{-1cm}
\ea

\begin{proposition}\label{prop: bundle str}
The induced bundle structure
(\ref{b-structure}) on $\p M$ and the Cauchy data set $C_0(D)$ determine a
 bundle $\tilde V$ on $M$ that is bundlemorphic to $V$. 
\end{proposition}

%Recall that the bundles $V$ and $\tilde V$ are isomorphic,
%if there is a fiber-preserving diffeomorphism $L:V\to \tilde V$
%that is a linear isomorphism at fibers.

\noindent {\em Proof:} Recall that by
Lemma \ref{lem: equivalence of data} we know $\Lambda_\Gamma$. By Lemma \ref{gauge3}
we can find, for any $y_0\in M$, a neighborhood $U$ 
and a system $\Phi$ of the eigenfunctions 
that define a local trivialization in $\pi^{-1}(U)$.
%
%defines a local trivialization of $V$
%compatible with the original structure of the vector 
%bundle $V$. Therefore, 
%As we have reconstructed the manifold $M$, it 
Thus it remains to prove that the boundary
data make it possible to evaluate 
the transition functions
\ba
T_{\Phi,U,\tilde \Phi,\tilde U}:(U\cap \tilde U)\times \C^d\to (U\cap \tilde 
U)\times \C^d
\ea
between any two  trivializations $A_{\Phi,U}$
and $A_{\tilde \Phi,\tilde U}$.

Let  $\hat f_+\in \F_+^\infty$, $\hat f_-\in \F_-^\infty$, and $t\in \R$.
By Lemma \ref{gauge3},
we can find decompositions
(\ref{27.11.2}) of ${\mathcal W}_+{\hat f_+}(t)$ and
${\mathcal W}_-{\hat f_-}(t)$ with respect to both
$A_{\Phi,U}$
and $A_{\tilde \Phi,\tilde U}$. Therefore, we
can find the representations of the function
\ba
w(y)={\mathcal W}_+{\hat f_+}(y,t) +{\mathcal W}_-{\hat f_-}(y,t),\quad
y\in U\cap \tilde U
\ea
in the both   trivializations $A_{\Phi,U}$
and $A_{\tilde \Phi,\tilde U}$.
According to Lemma \ref{C infinity}, the function $w$ 
runs over $C^\infty(M,V)$,  when   
 $\hat f_+$ runs over $\F_+^\infty$ and $\hat f_-$ over $ \F_-^\infty$. 
This makes it possible to compare representations
\ba
w(y) = 
\sum_{l=1}^d \alpha_l(y) \psi_{k(l)}(y)=
\sum_{l=1}^d \tilde \alpha_l(y) \psi_{\tilde k(l)}(y),\quad y \in U\cap \tilde U
\ea
when $w$ runs over $C^\infty(M,V)$ and, therefore, to determine
the transition function
$T_{\Phi,U,\tilde \Phi,\tilde U}$.
\proofbox\smallskip

\smallskip

\noindent {\bf  8.3.\ Construction of operators on $V$.} We are now in the position 
to complete the proof of the main Theorem \ref{main1}.

\noindent {\em Proof:} To this far, we have reconstructed $(M,g)$ upto
an isometry and a bundlemorphic copy $\tilde V$ of $V$.
To reconstruct the Hermitian and the Clifford module structures
and the operators $F$ and $D$  on $\tilde V$,
it is enough to find their representations 
in an arbitrary local trivialisation
$A_{\Phi,U}$. We will first consider the local  trivializations
$A_{\Phi,U}$ where
 $U\subset M^{\rm int}$  and 
 $\Phi$ consists of eigenfunctions of $D_{\Gamma,+}$ only.
 In this connection we will skip index $+$ in the following considerations.

We start with the reconstruction of the Hermitian 
structure. Let  $\hat f\in \F^\infty$,  
$y_0\in U$ and consider a ball $B(y_0,r)\subset U$. 
Denote by 
${\mathcal K}_r={\mathcal K}_r(\hat f)$ the set of those
 $\hat h\in \F^\infty$ for which
\beq\label{equals in U}
\bbra {\mathcal W}{\hat f}(T)-{\mathcal W}{\hat h}(T),
\lambda_y\ud_y\ccet=0,\ \hbox{for all }\lambda_y\in \pi^{-1}(y),\
 y\in B(y_0,r).\hspace{-2cm}
\eeq
Note that condition (\ref{equals in U}) is equivalent to the fact
that $ {\mathcal W}{\hat h}(T)|_{B(y_0,r)}=
{\mathcal W}{\hat f}(T)|_{B(y_0,r)}$.

Moreover, it follows from Lemma \ref{Lm2.15} that condition
 (\ref{equals in U}) can be directly verified from the boundary data.
 Therefore, using the boundary data we can evaluate
 \beq\label{104a}
 \| {\mathcal W}{\hat f}(T)\|_{L^2(B(y_0,r),V)}^2
=\inf_{\hh\in {\mathcal K}_r}
\|\hh\|_{\F}^2.
\eeq
As the manifold $(M,g)$ is already reconstructed
and ${\mathcal W}{\hat f}(T)\in C^\infty(M,V)$, 
the equation (\ref{104a}) makes it possible to find
\beq\label{Hermit}\quad\quad
| {\mathcal W}{\hat f}(y_0,T)|_{y_0}^2
=\lim_{r\to 0}
\frac {\|{\mathcal W}{\hat f}(T)\|_{L^2(B(y_0,r),V)}^2} 
{\hbox{vol}\,(B(y_0,r))}
\eeq
for an arbitrary $y_0\in U\subset M^{\rm int}$ and $\hat f\in \F^\infty_+$.
Again, as ${\mathcal W}{\hat f}(T)$ runs over ${\mathcal D}((D_{\Gamma,+})^\infty)$
when $\hat f$ runs over $\F^\infty_+$, if follows from
Lemma \ref{gauge2} that we can find 
the  Hermitian structure on $\pi^{-1}(y_0)$ in the trivialization $A_{\Phi,U}$.
As $y_0\in U$ is arbitrary,  we can determine the
Hermitian structure on $\tilde V|_{M^{\rm int}}$.

Next we consider projectors $\Pi_+:\pi^{-1}(y_0)\to \pi^{-1}(y_0)$,
$y_0\in M^{\rm int}$. Let
$\hat f\in \F^\infty$.
{\newwtextt By Theorem \ref{blacho} and Lemma \ref{gauge2} we first find, using the boundary data, 
$\hh \in {\mathcal F}$ such that $\Pi_+{\mathcal W}{\hat f}(T)={\mathcal W}{\hat h}(T)$.
As ${\mathcal W}{\hat f}(y_0,T)$ runs over $\pi^{-1}(y_0)$
when $\hat f$ runs over ${\mathcal F}$, this determines $\Pi_+$ in the trivialization $A_{\Phi,U}$. }
% Using
%Theorem \ref{blacho} and Lemma \ref{gauge2}
%we can find 
%the norms $|{\mathcal W}{\hat f}(y_0,T)|_{y_0}$ and $|\Pi_+{\mathcal W}{\hat f}(y_0,T)|_{y_0}$ analogously to (\ref{104a}).
%These determines the range of orthoprojector $\Pi_+$
%in the trivialization  $A_{\Phi,U}$, and therefore representation of $\Pi_+$
% in this trivialization.
Since $F=2\Pi_+-I$, the chirality operator $F$
is determined.

When $U\subset M^{\rm int}$ is small enough, there is 
an orthonormal frame $(e_1(y),\dots,e_n(y))$
on $TM|_{U}$. Let $y_0\in U$ and let $x(y)= (x^1(y),\dots, x^n(y))$, $x(y_0)=0$,
be  Riemannian normal coordinates
centered in $y_0$ with
  $e_l= \p_l$ at $y_0$.
By 
Theorem \ref{global control th} 
there are generalized
 sources ${\hat h}_{jl}\in \F^\infty,\, j=1,\dots, d,\, l=1,\dots,n$, such that
 \ba
u_{jl}(y,T):={\mathcal W}{\hat h}_{jl}(y,T)= x^l(y) \chi(y) \phi_{k(j)}(y),\quad y\in U.
\ea
Here $\chi\in C^\infty_0(U)$ is a cut-off function 
that is equal to $1$ 
in a neighborhood ${\tilde U}\subset U$
of $y_0$. 
Moreover, by  Lemma \ref{gauge2},
it is possible to find ${\hat h}_{jl}$ from the boundary data.

Representing the Dirac operator in this local trivialization (cf.\
(\ref{dirac 1}), (\ref{dirac-type, maximal})),
\ba
D u_{jl}(y,T) -\lambda_{k(j)}\,x^l(y)\, \phi_{k(j)}(y)= 
\sum_{p=1}^d \frac {\p}{\p x^p}(x^l)\, \gamma(e_p)\phi_{k(j)}(y),\quad y\in \tilde U.
\ea
In particular, at $y=y_0$
\beq
\label{27.11.3}
 D u_{jl}(y_0,T)= \gamma(e_l) \phi_{k(j)}(y_0).
\eeq
The left-hand side of (\ref{27.11.3}) coincides with
 $-i {\mathcal W}(\p_t\hat h_{jl})(y_0,T)$. By Lemma
 \ref{gauge2} the boundary data determines
 its representation in the local trivialization $A_{\Phi,U}$.
 Thus we can find $\alpha^{p}_{jl}(y_0)\in \C$ such that
 \ba
%\label{27.11.3b}
 \gamma(e_l) \phi_{k(j)}(y_0)=\sum_{p=1}^d 
 \alpha^{p}_{jl}(y_0)\,  \phi_{k(p)}(y_0).
\ea
This  provides
 the representation  of $\gamma(e_l)$ in the local trivialization 
 $A_{\Phi,U}$.
Thus we recover the Clifford module structure in $\tilde V|_{M^{\rm int}}$.

Next, we consider the representation of the Dirac-type operator
$D$. For any $\hat f\in \F^\infty$
we can find the representations of ${\mathcal W}{\hat f}(y,T)$
and $
D{\mathcal W}{\hat f}(y,T)=-i{\mathcal W}(\p_t \hat f)(T)$
in a local trivialization $A_{\Phi,U}$, $U\subset M^{\rm int}$. 
Recall that any  functions
$w\in C^\infty_0(U,\tilde V)$ can written as
 $w={\mathcal W}{\hat f}(T)$ with some $\hat f\in \F^\infty$.
 Therefore, the set of all pairs
\ba
\{(A_{\Phi,U}({\mathcal W}{\hat f}(T)|_U),A_{\Phi,U}(-i{\mathcal W}
(\p_t \hat f)(T)|_U)):\ \hat f\in \F^\infty\}
\ea
give the graph of 
the operator $D:C^\infty_0(U,\tilde V)\to C^\infty(U,\tilde V)$
in this local trivialization. This determines $D$ on $M^{\rm int}$.

To recover  the Hermitian and the Clifford module structures
and also the operators $F$ and  $D$  near $\p M$,
we recall that we have already found the bundle structure
of $\tilde V$ everywhere on $M$, i.e. the transition functions between local trivializations.
Consider  a local trivialization $A_{U,\Phi}$ 
where $U\cap \p M\not=\emptyset$ and  $\Phi=(\psi_{k(1)}(y),\dots,\psi_{k(d)}(y))$.
By the above considerations, we can find  
the representation  of the Hermitian and the Clifford module structures
and the operators $F$ and  $D$ on $U\cap M^{\rm int}$ in this trivialization.
However, these representations are smooth up to $\p M\cap U$,
and  we can continue them to $\p M\cap U$. Thus, we have
determined the  
Hermitian and the Clifford module structures
and the operators $F$ and  $D$ everywhere on $\tilde V$.
\hfill$\Box$

%Publish this later:
%\begin{remark}\rm 
%Similar to \cite{KLS1} we can show that Theorem \ref{main1} remains true 
%if, instead of
%$C_P$, or $\Lambda_{\Gamma}$, we are given  
%$\Lambda_{\Gamma}^T, \, T> 2 \rad(M)$. The response operator
%$\Lambda_{\Gamma}^T$ corresponds to finite time interval measurements,
%is  defined as
%\ba 
%\Lambda_{\Gamma}^T: f \mapsto u^f |_{\p M \times (0,T)},
%\quad f \in C^{\infty}_0(\p M \times (0,T), P_{\Gamma} W).
%\ea
%\end{remark}

\section{Inverse problem spectral problem.}
\noindent {\bf  9.1\ Time harmonic response operator.}
Let $\{\la_k^P, \,\phi_k^P|_{\p M}\}_{k=1}^{\infty}$ be the boundary spectral
 data corresponding
to the self-adjoint Dirac-type operator  (\ref{P}),
\bfo
D_Pu = Du, \quad
{\cal D}(D_P)= \{u \in H^1(M, V): \, P(u|_{\p M}) =0\},
\efo
where $P$ is a (possibly non-local) boundary condition satisfying 
conditions (\ref{Fredholm}) and 
(\ref{anticommutation}). Clearly, $\{\phi_k^P|_{\p M}\}_{k=1}^{\infty}$ span
the subspace $\,(I-P) H^{1/2}(\p M, W)$ and, therefore, determine $P$.

We intend to reconstruct, using the boundary spectral data for $D_P$, the fixed-frequency 
Cauchy data maps,
$\La^P_{\la}$,
\ba
%\label{Lambda}
\La^P_{\la}:\, P L^2(\p M, W) \rightarrow L^2(\p M, W), \quad
\La^P_{\la} f = u^f_{\la}|_{\p M}.
\ea
Here $u^f_{\la}$ is the solution to the boundary value problem,
\ba
Du^f_{\la} = \la u^f_{\la}, \quad P(u^f_{\la}|_{\p M})=f,
\quad \la \notin \sigma\,(D_P),
\ea
where $\sigma\,(D_P)$ is the spectrum of $D_P$.
$\La^P_{\la}$ is an integral operator with
its Schwartz kernel $\La^P_{\la}(x, y)$ formally 
defined as an expansion,
\bfo
\La^P_{\la}(x, y)=  \sum_{j=1}^{\infty} \frac{\phi_j^P(x) 
\,
\overline {\phi_j^P(y)}}
{\la-\la_j}.
\efo 
However, cf. \cite{NaSyU,KKLM} for the scalar Schr\"odinger operator, 
this series does not converge so that we
should regularize it. We use the same method as was used
in \cite{KKLM} for the Laplace operator, namely, we represent
the operator $\La^P_{\la}$ in the form, 
\ba
%\label{integral}
\La^P_{\la}f = \lim_{\mu \to +\infty}\La^P_{i \mu}f
+ \int^{\la}_{i \mu} \p_{\nu}(\La^P_{\nu}f)\, d\nu.
\ea
Observe that in the strong operator convergence 
in $H^{1/2}(\p M,\, W)$, 
\bfo
 \p_{\la}\La^P_{\la} = \lim_{A \to +\infty}  \p_{\la}\La^P_{\la,A},
\efo
where the Schwartz kernel of $\p_{\la}\La^P_{\la,A}$ is given by
\bfo
\p_{\la}\La^P_{\la,A}(x, y)= - \sum_{|\la_j^P| \leq A} 
\frac{\phi_j^P(x) \,
\overline {\phi_j^P(y)}}
{(\la-\la_j)^2}.
\efo
 Thus it 
remains to find $ \lim_{\mu \to +\infty}\, \La^P_{i \mu}$.
\begin{lemma}
\label{limit}
Let $D_P$ be an operator of form (\ref{P}) where $P$ satisfies 
conditions
(\ref{anticommutation}),  (\ref{Fredholm}). Then,
\beq
\label{WKB}
 \lim_{\mu \to +\infty}\, \La^P_{i \mu} \,f =
f+i \gamma(N)f.
\eeq
\end{lemma}
\noindent {\em Proof:}
Consider the boundary value problem 
\beq
\label{Dirichlet}
D\, u^f_{i\mu}= i\mu\, u^f_{i\mu}, \quad
P (u^f_{i\mu}|_{\p M})= f \in P C^{\infty}(\p M,\, W).
\eeq
We introduce the boundary geodesic coordinates $x(y)=(x'(y), x_n(y))$
where $x_n= \hbox{dist}(y, \, \p M)$ and $x'$ are the coordinates
of the nearest boundary point.
Using a special trivialisation of $V$ near $\p M$,
see Appendix, formulae (\ref{newform}) and (\ref{A}), 
 the Dirac-type operator $D$
takes the form
\ba
%\label{in bnc}
& &   \quad\ \  D = \gamma(N) \left(\p_n + {\it A}(x',x_n, D')\right), \\
\nonumber
& & \hspace{-12mm}{\it A}(x',x_n, D')= -\gamma(N)\sum_{\alpha=1}^{n-1} \gamma(e_\a) \triangledown_{\a}
+ \frac{n-1}{2}  H(x', x_n) - \gamma(N) Q(x', x_n) .\hspace{-5mm}
\ea
Here $ H(x', \tau)$ is the mean curvature of the surface
$\p M_{\tau} = \{y\in M: \dist(x,\p M)=\tau\}$ at the point $y$ with coordinates
$(x',\tau)$.

We look for the solution of (\ref{Dirichlet}) near 
$\p M$, using  a formal
WKB expansion,
\beq
\label{WKB1}
u^f_{i\mu}(x',x_n) \sim e^{-\mu x_n} \sum_{k=0}^{\infty}
\frac{1}{\mu^k}\, u_k(x', x_n).
\eeq
Substitution of (\ref{WKB1}) to (\ref{Dirichlet}) gives rize to
a recurrent system of equations with respect to the powers of $\mu$,
\beq
\label{recurrent}
-(i + \gamma(N))\, u_k+{\it A}(x',x_n, D') u_{k-1}
+\gamma(N) \p_nu_{k-1}
=0,
\eeq
where $u_{-1}=0$. This system has to be looked at together with the boundary conditions,
\beq
\label{boundarycond}
P (u_k|_{x_n=0}) = \left\{ \begin{array}{l} f,\quad \hbox{for }k=0,\\
0,\quad \hbox{for }k\not=0.\end{array}\right.  
\eeq
Equation (\ref{recurrent}) tells us that,
for any $x$,  $u_0$ is the eigenspace of $\gamma(N)$ corresponding to
the eigenvalue $-i$. We denote the eigenprojectors  
of $\gamma(N)$ corresponding to
the eigenvalues $-i$ and $+i$ by 
$p^{-}_{\gamma}$ and $p^{+}_{\gamma}$,
\ba%\label{eigen pr.}
p^{-}_{\gamma} =\frac 12 (1 +i \gamma(N)),\quad
p^{+}_{\gamma} =\frac 12 (1 -i \gamma(N)).
\ea
Applying projector $P$ to (\ref{recurrent}) with $k=0$ and using  
equation
 $\gamma(N) P= (I-P) \gamma(N)$, we get 
\bfo
i P u_0+ \gamma(N) (I-P) u_0 =0.
\efo
Letting $x_n=0$ in the above equation and using (\ref{boundarycond}),
we obtain equation (\ref{WKB}).
To justify it, we, however,  need to consider further terms in 
the WKB-expansion
(\ref{WKB1}).

If we apply the projector $p^{-}_{\gamma}$ to (\ref{recurrent})
with $k=1$,
we see that
\beq
\label{0-order}
-i \p_n p^{-}_{\gamma} u_0+p^{-}_{\gamma}\,{\it A}(x',x_n, D')\,
p^{-}_{\gamma}u_0=0,
\eeq
where we use that $\p_n \gamma(N)=0$ and $u_0=p^{-}_{\gamma}u_0$.
Since $\gamma (e_{\a})\,\gamma (N)+\gamma (N)\,\gamma (e_{\a})=0$,
the differential operator $p^{-}_{\gamma}\,{\it A}(x',x_n, D')\,
p^{-}_{\gamma}$  is of the order $0$. Therefore, (\ref{0-order}) is an 
ordinary differential equation 
for $p^{-}_{\gamma}u_0$  along the normal geodesic with the
initial condition  given by (\ref{WKB}). This determines $u_0$.
Next
we apply $p^{+}_{\gamma}$ to (\ref{recurrent}) with $k=1$ to obtain
$p^{+}_{\gamma} u_1$,
\bfo
2i p^{+}_{\gamma} u_1 =p^{+}_{\gamma}\,{\it A}(x',x_n, D')u_0.
\efo
Together with (\ref{boundarycond}) with $k=1$, this determines
$u_1|_{x_n=0}$. Further considerations follow the trend above: applying 
$p^{-}_{\gamma}$ to (\ref{recurrent}) with $k=2$, we obtain an ordinary 
 differential equation along the normal geodesic for $p^{-}_{\gamma} u_1$;
applying $p^{+}_{\gamma}$ to (\ref{recurrent}) with $k=2$ we determine
$p^{+}_{\gamma} u_2$ and, by (\ref{boundarycond}) with $k=2$, the value
$u_2|_{x_n=0}$, etc.

Take a finite sum $\exp(-\mu x_n) \sum_{0}^l \mu^{-k}u_k$
in  (\ref{WKB1}) and  
cut it by multiplying with
a smooth function $\chi(x_n)$ equal to $1$ near $\p M$.  
Since
$
P(\sum_{k=0}^l \mu^{-k} u_k|_{\p M})=f,
$
and 
$
\hbox{dist}\left(i\mu,\, \sigma(D_P)\right) \geq |\mu|,
$
we see 
that
\bfo
\|u^f_{i\mu}-e^{-\mu x_n} \sum_{k=0}^l \mu^{-k} u_k \chi\|_{H^1(M)} 
\leq c_l \mu^{-l-1}.
\efo
implying the desired result on $\p M$.
\hfill$\Box$

%Having in our disposal $\La^P_{\la}$ and, therefore, the Cauchy data
%for the solutions to $Du^f_{\la} = \la u^f_{\la}$, it is possible, 
%using the
%knowledge of $\gamma(N)|_{\p M}, \, F|_{\p M}$, to find the Cauchy data map
%$\La^L_{\la}$ which corresponds to the local boundary condition,
%\bfo
%P^L u= \frac12 (I+ \gamma(N) \circ F) u = 0.
%\efo
%Using the same considerations as in \cite{KKLM} or \cite{KKL}, Sec. 4.1,
%we can find from $\La^L_{\la}, \, \la \notin \hbox{spec}( D^L)$ the boundary spectral data
%$\{\la_k,\, \phi_k^L|_{\p M}\}$ for the Dirac operator with the 
%local 
%boundary
%condition.

\smallskip

\noindent {\bf  9.2.\ From time harmonic data to time domain data.}
Having $\Lambda_{\la}^P$ in our disposal, we reconstruct the 
response operator $\Lambda_P$, see (\ref{response}). 
%Indeed, as 
%$D_P$ is self-adjoint, the initial boundary value problem
%(\ref{IBVP0})-(\ref{IBVP0line 2}) is well-defined. 
To this end we observe that the set
\bfo
\{u^f|_{\p M \times \R_+}: \, u^f \,\, \hbox{solves (\ref{IBVP0})
with some}\,\, f \in P\Cnull(\p M \times \R_+,\, W)\},
\efo 
coincides with the Cauchy data set $C_0(D)$. Due to causality,
to reconstruct $\Lambda_P$ it is sufficient to find its values
only on $PC^{\infty}_0(\p M \times \R_+,\, W)$.

Take $f \in PC^{\infty}_0(\p M \times \R_+,\,  W)$, and let
$T^f\in \R_+$ be such  that  $f=0$ for $t> T^f$. Therefore, 
$u^f(t) \in {\mathcal D}(D_P)$ for $t >T^f$. Moreover, the graph norm 
of the wave $u^f(t)$,
\bfo
\|u^f(t)\|^2_{{\mathcal D}(D)}=
\|Du^f(t)\|^2_{L^2(M, V)} + \|u^f(t)\|^2_{L^2(M, V)},
\efo
is independent of $t, \, t> T^f$. Since
this norm is equivalent to the 
$H^1(M, V)-$norm  on ${\mathcal D}(D_P)$,  we see that 
\beq
\label{estimate}
u^f \in C_{0, b}(\R_+, \, H^{1}(M, V)), \
u^f|_{\p M\times \R_+} \in C_{0, b}(\R_+, \, H^{1/2}(\p M, W)),\hspace{-2cm}
\eeq
where  $C_{0, b}(\R_+,\, X)$ is the class of bounded  
functions with values in $X$, which are   equal to $0$
near $t=0$.

Consider the Fourier transform,
\bfo
{\widehat u}^f(x, \la) = \int_\R \exp(-i \la t) u^f(x, t) dt.
\efo
By (\ref{estimate}), ${\widehat u}^f(\cdot, \la)$ is an analytic
$H^1(M, V)$-valued  function
of $\la$ in the lower half plane, $\la = \nu +i \mu, \, \mu <0$. 
For each $\la,\,$ ${\widehat u}^f(\cdot, \la)$
is the solution of the Dirichlet problem (\ref{Dirichlet}) with 
$\la$ instead of $i \mu$ and
${\widehat f}(\la)$ instead of $f$ in the boundary condition.
From  $\Lambda_{\la}^P$ we find ${\widehat u}^f(\la)|_{\p M}$.
It follows from (\ref{estimate}) that ${\widehat u}^f(\la)|_{\p M}$
is bounded along the line $\{\nu +i \mu: \, \nu \in \R\}$
 for any $\mu <0$.
Thus,
\bfo
u^f|_{\p M \times \R_+} =
\frac{1}{2 \pi} \exp(-\mu t)\int_{\R} \exp(i \nu t) 
\Lambda^P_{\nu+i\mu}{\widehat f}(\nu +i \mu)|_{\p M} d\nu.
\efo

Summarizing, we see that the induced bundle structure on $\p M$
and the boundary spectral data 
$\{\la_k, \,\phi_k|_{\p M}\}_{k=1}^{\infty}$ of a
Dirac-type operator $D_P$, with $P$ satisfying (\ref{Fredholm}),
(\ref{anticommutation}), determine the Cauchy data set, $C_0(D)$.
By Theorem \ref{main1} this implies that the boundary spectral data determine
the Riemannian manifold,
$(M, g)$, up to an isometry, and the vector bundle, $V$, with its Hermitian
and Clifford module structures, as well as the chirality operator, $F$ and
Dirac-type operator, $D$, up to a Dirac bundlemorphism. 
Finally, since $\phi_k|_{\p M},$ $k=1,2,\dots$, span the subspace
$P\,C^{\infty}(\p M,W)$,
we  can  determine the projector $P$. Thus we can determine $D_P$.
This proves Theorem \ref{main2}.\hfill$\Box$

\section{Dirac-type operator in $\R^3$}
In this section we apply our results to the classical example
of the canonical Dirac-type operator in a bounded domain $M \subset \R^3$.

Let $V=M \times \C^4$ be the trivial bundle  over $M$, each fiber endowed
with the standard Hermitian inner product of $\C^4$.
The  unperturbed Euclidean Dirac operator, $D_0$ has the form
\beq
\label{standard}
D_0u =i \sum_{k=1}^3 \alpha_k \p_k u + mc^2\alpha_0 u, \quad
u \in H^1(\R^3; \C^4).
\eeq 
Here $\alpha_{\nu},\, \nu  =0, \dots,3,$ are the standard Dirac matrices,
\ba
\alpha_0=
\left(
\begin{array}{cc}
I & 0\\
0 & -I
\end{array}
\right),
\quad \alpha_k =
\left(
\begin{array}{cc}
0 & \sigma_k\\
\sigma_k & 0
\end{array}
\right),\quad k=1,2,3,
\ea
where $\sigma_k$ are the Pauli $2\times 2$ matrices.
%\HOX{What do Pauli matrices satisfy, reference}
Then
\ba
%\label{standard anticom}
 \alpha_{\nu} \alpha_{\mu}+\alpha_{\mu} \alpha_{\nu}  =
\left\{ \begin{array}{l} 2I ,\quad \hbox{for }\nu=\mu,\\
0,\quad \hbox{for }\nu\not=\mu.\end{array}\right. 
, \quad \nu, \mu=0,\dots, 3.
\ea
and the operator (\ref{standard})
  falls into the category of the Dirac-type operators
considered in this paper with
$
F =  \alpha_0 \alpha_1  \alpha_2  \alpha_3$ and 
$\gamma(e_k) = i \alpha_k.$ 

A perturbation $Q$ of a physical nature,
consisting of a zero components of the Lorentz 
scalar and vector potentials, see e.g.\  \cite{NWvD},
defines a perturbed Dirac operator  $D_{{\bf a},q}=D_0+Q$. Namely,
\beq\label{Dirac, euclidean}
D_{{\bf a},q} u =i \sum_{k=1}^3 \alpha_k (\p_k+ia_k(x)) u + 
mc^2  \alpha_0 u + \alpha_0 q(x) u,
\eeq
where $a_k(x),q(x)$ are $C^\infty$-smooth real-valued functions
which correspond to the magnetic vector potential
${\bf a}(x)=a_1dx^1+a_2dx^2+a_3dx^3\in \Omega^1(M)$
 and the scalar potential $q(x)$.
Then the potential $Q$ satisfies conditions 
(\ref{Q self-adjoint}) and (\ref{respect}) and, applying
the previous constructions of the paper,
we obtain the following theorem:

\begin{theorem}
\label{main Euclidean}
Let $M\subset \R^3$ be a bounded domain with a smooth boundary.
Let $D_{{\bf a},q}$ and $D_{\tilde{\bf a},\tilde q}$
 be two Dirac-type operators of 
 form (\ref{Dirac, euclidean}) in $M$. Assume that the
Cauchy data sets of these operators are the same, i.e.
\ba
C_0(D_{{\bf a},q})=C_0(D_{\tilde{\bf a},\tilde q}).
\ea
Then $\tilde q= q$ and
\beq\label{tilde A and A}
\tilde {\bf a}-{\bf a} = d \Phi + \sum_{m=1}^{b_1} n_m h_m,
\quad n_m \in \Z,
\eeq
where $\{h_m\}_{m=1}^{b_1}$ is a basis of the relative  harmonic $1-$forms on $M$,
$b_1$ is the first relative Betti number of $M$ and $d\Phi$ 
is the exterior differential of a function $\Phi
\in C^\infty(M,\R)$,  $\Phi|_{\p M}=0$.
\end{theorem}
 
\noindent
The 1-forms $h_m$ are normalized so
that if $\p M$ has the external component $\Sigma_0$ and
the internal components $\Sigma_j$, $j=1,\dots,b_1$, then
for a path $\eta=\eta([0,1])\subset M$, $\eta(0)\in \Sigma_0$,
and $\eta(1)\in \Sigma_j$ we have $\int_\eta h_m=\delta_{jm}.$
\smallskip

\noindent {\em Proof:} 
{\newwtext Clearly, the standard Hermitian structure of $\C^4$, 
$\gamma(e_k) = i \alpha_k$, and the canonical differential, $\triangledown_{e_k}= \p_k$. 
define the structure of a Dirac bundle in $M \times \C^4$ with
$D_{{\bf a},q}$ and $
D_{\tilde{\bf a},\tilde q}$ being  Dirac-type operators on this bundle. }
%We consider the operators $D=D_{{\bf a},q}$ and $\tilde D=
%D_{\tilde{\bf a},\tilde q}$ as Dirac operators defined
%on the same bundle $M\times \C^4$ that is endowed with
%the standard Hermitian structure of $\C^4$,
%the chirality operator $
%F =  \alpha_0 \alpha_1  \alpha_2  \alpha_3$ and 
%the action $\gamma(e_k) = i \alpha_k$  of the Clifford algebra.

By Theorem \ref{main1}, there is a bundlemorphism,
$L:M\times \C^4\to M\times \C^4$ that is a fiberwise
isomorphism and, therefore, is given by a smooth matrix-valued function $L(x)$ which
 is invertible for all $x\in M$. 
Applying again Theorem \ref{main1}, we see that
 $D=L^{-1}\tilde D L$, $F= L^{-1} \tilde F L$,  and 
 $\gamma(e_k )=L^{-1}\gamma(e_k )L$ for $k =1,2,3$. These imply that
 $\a_{\nu} L =L \a_{\nu},\, \nu =0, \dots, 3,$ which means, in turn, that 
 $L(x)$ is a diagonal matrix.
As the bundlemorphism $L$ preserves the Hermitian structure 
in each
fiber, there is a complex valued function $\kappa(x)$
such that
\beq\label{kappa-morphism}
L(x)= \kappa(x) I, \quad |\kappa(x)|=1.
\eeq
 Since the induced bundle structures on the boundary coincide,
$L(x)=I$ for $x\in \p M$. Together with (\ref{kappa-morphism}) this yields that,
 cf.  \cite{Ku3},
$
\kappa(x)= \exp(i\Psi(x)),
$
with $\Psi\in C^\infty(M,\R)$ such that
$\Psi(x)=0$ on the external component of $\p M$ and
$\Psi(x)=2\pi n_m$, $n_m\in \Z$ on the internal $m$-th component of $\p M$.
The family of bundlemorphisms (\ref{kappa-morphism}) 
are precisely
those bundlemorphisms that preserve the bundle
structure and, in particular,
  the structure (\ref{Dirac, euclidean})
of the Dirac operator. 

Since  $D=L^{-1}\tilde D L$, where $L$ is of  form
(\ref{kappa-morphism}), we see that
$\tilde q=q$ and
\beq\label{gauge to a}
\tilde {\bf a}={\bf a} + d\Psi.
\eeq
We obtain formula (\ref{tilde A and A})
from  (\ref{gauge to a}) by writing $d\Psi=d\Phi
 + \sum_{m=1}^{b_1} n_m h_m$ where $\Phi|_{\p M}=0$,
$b_1$ is the first relative Betti number of $M$ and
$h_m$ is a basis of the   harmonic $1-$forms on $M$
satisfying the relative boundary condition (\ref{relative}).
\proofbox\smallskip

Theorem \ref{main Euclidean} 
 is a improvement, for the considered class of 
perturbations and domain $M\subset \R^3$, of the results in  \cite{Nakamura}.
Note that Theorem \ref{main Euclidean} describes, for the Dirac operator, 
the relation between the inverse problems and the Aharonov-Bohm effect.
(For the previous results 
on the Aharonov-Bohm effect for a Schr\"odinger operator,
see e.g.\  \cite{Es,We} for the Euclidean Laplacian and \cite{KK2}
for a more general scalar operators on a Riemannian manifold.) 
In physical
terms,  Aharonov-Bohm effect means that
not only the 
the magnetic field $d{\bf a}$, i.e.
$\hbox{curl}({\bf a})$, but also the 
magnetic vector potential ${\bf a}(x)$ has an effect on the
 boundary measurements. Theorem \ref{main Euclidean} characterizes
what part of ${\bf a}$ besides $d{\bf a}$ can be observed from
external measurements.

\subsection*{Appendix. Local Boundary Condition}
In this appendix we provide a trivialisation 
of $V$ near the boundary needed in Sections 2 and 9, and
outline the proof that the Dirac-type 
operator with the local boundary condition is self-adjoint. 
Here, we consider the unperturbed Dirac operator
\ba
%\label{localbc}\quad\quad
D_{\Gamma}u = D_0u, \quad {\mathcal D}(D_{\Gamma}) = \{u \in H^1(M, V): \, 
P_{\Gamma}(u|_{\p M})=0\},
\ea
where $D_0$ is defined by (\ref{dirac 1}) and
the projector $P_{\Gamma}$ by (\ref{Gamma}), 
(\ref{PGamma}).

First we consider the needed trivialisation of $V$.
In the following, we do not change notations for 
$\gamma(v), v \in T_x(M), \, F,$ etc.\ when changing from one
trivialisation to another.

Let $U \subset M$ be a coordinate chart near $\p M$ with boundary
normal coordinates, $x(y)=(x'(y), x_n(y)).$ 
We may assume 
that $x(U)=B'\times[0,\e]$, where $B'\subset \R^{n-1}$ is open.
Let $w_j(x)\in \pi^{-1}(x)$, $x\in U\cap \p M$, be smooth sections 
such that $\{w_j(x)\}_{j=1}^d$ is an orthonormal basis of  $\pi^{-1}(x)$,
$x\in U\cap \p M$.
We can continue $w_j$ along the normal geodesics so that
$\nabla_N w_j=0$, and obtain sections $w_j \in C^{\infty}U, V)$
that form an orthonormal basis for $\pi^{-1}(x)$,
$x\in U$. The map \ba
\Psi:\bigg(x,\sum_{j=1}^d a_j(x)w_j(x)\bigg)\mapsto (x,
a_1(x),\dots, a_d(x))
\ea 
defines a local trivialisation
$\Psi:\pi^{-1}(U) \to  U \times \C^d$.
In this trivialisation $\triangledown_N u$ is represented by
$\p_n u$.
% and we have
%$
%\p_n \gamma(N)=0.
%$

Denoting $\Psi(x,v)=(x,\psi(x,v))$, we can define another
trivialisation $\Phi(x,v)=(x,\omega(x)\psi(x,v))$,
where in the boundary normal coordinates
\bfo
\omega(x', x_n)= 
\exp\left(-\frac{n-1}{2} \int_0^{x_n}H(x', \tau) d \tau\right)
\efo
with $H(x',\tau)$ being the mean curvature of the hypersurface
$\p M_{\tau}= \{x\in M: \dist(x,\p M)=\tau\}$. In this
trivialisation  $D_0$ has the form
\beq
\label{newform} 
D_0=\gamma(N)[\p_n + A(x', \p_{x'}; x_n)] ,
\eeq
where $A(x_n)= A(x', \p_{x'}; x_n)$ is the first-order differential
operator with respect to $x'$ which depends on $x_n$ as a parameter.
Namely,
\beq
\label{A}
A(x', \p_{x'}; x_n) = 
- \gamma(N) \sum_{\a=1}^{n-1} \gamma(e_{\a})\triangledown_{\a} +
\frac{n-1}{2} H(x',x_n).
\eeq
Here, $\triangledown_{\a}=\triangledown_{e_\a}.$
%An important property of the trivialisation $\Phi$ is that
%$D_0$ is represented as
%a formally self-adjoint operator (\ref{newform}) in
%$L^2([0,\e],dx_n;H)$, where $H=L^2(B',\C^d)$ is a Hilbert space
%that does not depend on $x_n$ and $dx_n$ is the 
%Lebesgue measure on interval $[0,\e]$.
The obtained operator $A(x',\p_{x'};0)$ is called the {\it hypersurface
Dirac operator.}

% $H=L^2(\p M,V|,d\mu)$ where the measure 
% the covariant derivative 
%
%, we can consider
%(\ref{newform}) as a formally self-adjoint operator on
%$L^2([0,\e],dx_n;H)$, $H=L^2(\p M,V|,d\mu)$ where the measure 
%$d\mu$ is independent of $x_n$ and is given by
%$d\mu=\det(g_{\alpha\beta}(x',0))^{1/2}dx_1'\dots dx_{n-1}'$,
%where $g_{\alpha\beta}(x',0)$ is the metric tensor of $\p M$.

Next we show  that
$D_{\Gamma}$ fits into the framework of the theory developed in  
\cite{BrL} implying its self-adjointness and regularity. (For other proofs of
the self-adjointness see \cite{FS,HSZ1}.)
To this end, we use the  above
 trivialisation $\Phi$ of $V$ near $\p M$ so that
$D_0$ has form (\ref{newform}), (\ref{A}).
We first show  that $A=A(\tau):L^2(H)\to L^2(H)$ satisfies
\beq
\label{cond1}
\gamma(N) A+ A\gamma(N)=0 \quad\hbox{and}\quad
A(0)=A^*(0).
\eeq

To prove the first identity in (\ref{cond1}), we introduce the Riemannian 
normal coordinates 
 on $\p M_{\tau}$ centered at $x'_0$ which are associated with the 
orthonormal frame
$(e_1, \dots, e_{n-1})$. Then,
\beq
\label{computation}\quad\quad
\gamma(N) A =-\gamma(N) \sum_{\a=1}^{n-1} \gamma(N) \gamma(e_{\a})
\triangledown_{\a} +
\frac{n-1}{2} H(x',x_n)\gamma(N).
\eeq
At the point $(x'_0, \tau)$,
\ba
& &\triangledown_{\a}(\gamma(N)u)= \gamma(N)\triangledown_{\a} u
+ \gamma(\triangledown_{\a}N) u
= \gamma(N)\triangledown_{\a} u -
\sum_{\beta=1}^{n-1}S_{\a\beta} \gamma(e_{\beta}) u,
\ea
where $S_{\a\beta}$ is the 2nd fundamental form of $\p M_{\tau}$ in 
the introduced normal coordinates centered at $x'_0$. Returning to 
(\ref{computation}), this implies that
\bfo
\gamma(N) A = 
\gamma(N)\left(
  \sum_{\a, \beta=1}^{n-1}S_{\a\beta}\gamma(e_{\a}) \gamma(e_{\beta})
+ \sum_{\a=1}^{n-1}  \gamma(e_{\a}) \triangledown_{\a} \gamma(N)
+ \frac{n-1}{2} H\right).
\efo
By (\ref{clifford}), 
$\gamma(e_{\a}) \gamma(e_{\beta})+\gamma(e_{\beta}) \gamma(e_{\a})=-2\delta_{\a \beta}$
and $S_{\a\beta}$ is symmetric. Therefore,
the above equation and (\ref{computation}) yield 
%that
%\bfo
%\gamma(N) A \hspace{-1mm}=\hspace{-1mm} 
%\gamma(N)\left( 
%-(n-1)H+\sum_{\a=1}^{n-1}  \gamma(e_{\a}) \triangledown_{\a} \gamma(N)
%+ \frac{n-1}{2} H\right) \hspace{-1mm}=\hspace{-1mm}
%- A  \gamma(N),
%\efo
%providing 
the first identity in (\ref{cond1}).

Similar arguments, employing that in the boundary
 normal coordinates, 
$
%beq
%label{second}
\triangledown_{\a} e_{\beta}= S_{\a\beta} N,$ $N=\p_n,
%eeq
$
show the second identity in (\ref{cond1}). 
%Indeed,
%\bfo
%& &A^*(0)= \sum_{\a=1}^{n-1}\triangledown_{\a}  \gamma(e_{\a})\gamma(N)
% - \frac{n-1}{2} H 
%\\
%\nonumber
%& &=
% \sum_{\a=1}^{n-1} \gamma(e_{\a}) \triangledown_{\a} \gamma(N)
%+(n-1) H  - \frac{n-1}{2} H 
%\\
%\nonumber
%& &
%=
%- \gamma(N) \sum_{\a=1}^{n-1} \gamma(e_{\a}) \triangledown_{\a}
%-(n-1) H +(n-1) H  - \frac{n-1}{2} H 
%= A(0).
%\efo

It  remains to show that $P_{\Gamma}$ satisfies
(\ref{anticommutation}) and  (\ref{Fredholm}).
As
$
P_{\Gamma} = \frac12 (I+\Gamma),
$
it follows from
(\ref{clifford}), (\ref{chirality2}) that $P_{\Gamma}$ satisfies
(\ref{anticommutation}). To check (\ref{Fredholm}), let $P_{APS}$
be a spectral projector of $A(0)$ satisfying
(\ref{APS}). First equation
of (\ref{cond1}) and (\ref{chirality2})
imply that $A(I+\Gamma)=(I-\Gamma)A.$
Thus if  $\la_k$ is an eigenvalue of $A(0)$ with an 
eigenfunction $\phi_k$, then
\bfo
\psi_k = P_{\Gamma} \phi_k -(I-  P_{\Gamma}) \phi_k
\efo
is an eigenfunction for the eigenvalue $- \la_k$. 
In particular, when $ \la_k \not= 0$,
$\bra\bra\phi_k, \psi_k \cet\cet=0,$ so that
$
\|P_{\Gamma} \phi_k\|= \|(I-P_{\Gamma}) \phi_k\|$ and $
\|\phi_k\|= \|\psi_k\|.
$
Consider the decomposition,
$
L^2(\p M, W) = {\Bbb H}_1 \oplus {\Bbb H}_2,
$
where ${\Bbb H}_2$ is the  $0-$eigenspace of
$A(0)$. The above equations imply that ${\Bbb H}_1, \,{\Bbb H}_2$
are invariant subspaces for $P_{\Gamma}$ and
\bfo
\|P_{\Gamma} - P_{APS}\| = \frac{1}{{\sqrt 2}} \quad \hbox{on} \,\,
{\Bbb H}_1.
\efo
As ${\Bbb H}_2$ is finite-dimensional, then
$1$ and $-1$ are not in the essential spectrum of $P_{\Gamma} -  P_{APS}$.
By \cite{BrL}, Proposition 3.3, $(P_{\Gamma}, \, P_{APS})$ is a 
Fredholm pair, so that $D_{\Gamma}$ is self-adjoint.

\end{document}